\def\fmtname{AmS-TeX}

\def\fmtversion{2.1}
\catcode`\@=11
\ifx\amstexloaded@\relax\catcode`\@=\active
  \endinput\else\let\amstexloaded@\relax\fi
\newlinechar=`\^^J
\def\W@{\immediate\write\sixt@@n}
\def\CR@{\W@{^^J\fmtname - Version \fmtversion^^J%
COPYRIGHT 1985, 1990, 1991 - AMERICAN MATHEMATICAL SOCIETY^^J%
Use of this macro package is not restricted provided^^J%
each use is acknowledged upon publication.^^J}}
\CR@ \everyjob{\CR@}
\message{Loading definitions for}
\message{misc utility macros,}
\toksdef\toks@@=2
\long\def\rightappend@#1\to#2{\toks@{\\{#1}}\toks@@
 =\expandafter{#2}\xdef#2{\the\toks@@\the\toks@}\toks@{}\toks@@{}}
\def\alloclist@{}
\newif\ifalloc@
\def\showallocations{{\def\\{\immediate\write\m@ne}\alloclist@}\alloc@true}
\def\alloc@#1#2#3#4#5{\global\advance\count1#1by\@ne
 \ch@ck#1#4#2\allocationnumber=\count1#1
 \global#3#5=\allocationnumber
 \edef\next@{\string#5=\string#2\the\allocationnumber}%
 \expandafter\rightappend@\next@\to\alloclist@}
\newcount\count@@
\newcount\count@@@
\def\FN@{\futurelet\next}
\def\DN@{\def\next@}
\def\DNii@{\def\nextii@}
\def\RIfM@{\relax\ifmmode}
\def\RIfMIfI@{\relax\ifmmode\ifinner}
\def\setboxz@h{\setbox\z@\hbox}
\def\wdz@{\wd\z@}
\def\boxz@{\box\z@}
\def\setbox@ne{\setbox\@ne}
\def\wd@ne{\wd\@ne}
\def\iterate{\body\expandafter\iterate\else\fi}
\def\err@#1{\errmessage{AmS-TeX error: #1}}
\newhelp\defaulthelp@{Sorry, I already gave what help I could...^^J
Maybe you should try asking a human?^^J
An error might have occurred before I noticed any problems.^^J
``If all else fails, read the instructions.''}
\def\Err@{\errhelp\defaulthelp@\err@}
\def\eat@#1{}
\def\in@#1#2{\def\in@@##1#1##2##3\in@@{\ifx\in@##2\in@false\else\in@true\fi}%
 \in@@#2#1\in@\in@@}
\newif\ifin@
\def\space@.{\futurelet\space@\relax}
\space@. %
\newhelp\athelp@
{Only certain combinations beginning with @ make sense to me.^^J
Perhaps you wanted \string\@\space for a printed @?^^J
I've ignored the character or group after @.}
{\catcode`\~=\active 
 \lccode`\~=`\@ \lowercase{\gdef~{\FN@\at@}}}
\def\at@{\let\next@\at@@
 \ifcat\noexpand\next a\else\ifcat\noexpand\next0\else
 \ifcat\noexpand\next\relax\else
   \let\next\at@@@\fi\fi\fi
 \next@}
\def\at@@#1{\expandafter
 \ifx\csname\space @\string#1\endcsname\relax
  \expandafter\at@@@ \else
  \csname\space @\string#1\expandafter\endcsname\fi}
\def\at@@@#1{\errhelp\athelp@ \err@{\Invalid@@ @}}
\def\atdef@#1{\expandafter\def\csname\space @\string#1\endcsname}
\newhelp\defahelp@{If you typed \string\define\space cs instead of
\string\define\string\cs\space^^J
I've substituted an inaccessible control sequence so that your^^J
definition will be completed without mixing me up too badly.^^J
If you typed \string\define{\string\cs} the inaccessible control sequence^^J
was defined to be \string\cs, and the rest of your^^J
definition appears as input.}
\newhelp\defbhelp@{I've ignored your definition, because it might^^J
conflict with other uses that are important to me.}
\def\define{\FN@\define@}
\def\define@{\ifcat\noexpand\next\relax
 \expandafter\define@@\else\errhelp\defahelp@                               
 \err@{\string\define\space must be followed by a control
 sequence}\expandafter\def\expandafter\nextii@\fi}                          
\def\undefined@@@@@@@@@@{}
\def\preloaded@@@@@@@@@@{}
\def\next@@@@@@@@@@{}
\def\define@@#1{\ifx#1\relax\errhelp\defbhelp@                              
 \err@{\string#1\space is already defined}\DN@{\DNii@}\else
 \expandafter\ifx\csname\expandafter\eat@\string                            
 #1@@@@@@@@@@\endcsname\undefined@@@@@@@@@@\errhelp\defbhelp@
 \err@{\string#1\space can't be defined}\DN@{\DNii@}\else
 \expandafter\ifx\csname\expandafter\eat@\string#1\endcsname\relax          
 \global\let#1\undefined\DN@{\def#1}\else\errhelp\defbhelp@
 \err@{\string#1\space is already defined}\DN@{\DNii@}\fi
 \fi\fi\next@}

\def\predefine#1#2{\let#1#2}
\def\undefine#1{\let#1\undefined}
\message{page layout,}
\newdimen\captionwidth@
\captionwidth@\hsize
\advance\captionwidth@-1.5in
\def\pagewidth#1{\hsize#1\relax
 \captionwidth@\hsize\advance\captionwidth@-1.5in}

\def\hcorrection#1{\advance\hoffset#1\relax}
\def\vcorrection#1{\advance\voffset#1\relax}
\message{accents/punctuation,}

\let\graveaccent\`
\let\acuteaccent\'
\let\tildeaccent\~
\let\hataccent\^
\let\underscore\_
\let\B\=
\let\D\.
\let\ic@\/
\def\/{\unskip\ic@}
\def\textfonti{\the\textfont\@ne}
\def\t#1#2{{\edef\next@{\the\font}\textfonti\accent"7F \next@#1#2}}
\def~{\unskip\nobreak\ \ignorespaces}
\def\.{.\spacefactor\@m}
\atdef@;{\leavevmode\null;}
\atdef@:{\leavevmode\null:}
\atdef@?{\leavevmode\null?}
\edef\@{\string @}
\def\relaxnext@{\let\next\relax}
\atdef@-{\relaxnext@\leavevmode
 \DN@{\ifx\next-\DN@-{\FN@\nextii@}\else
  \DN@{\leavevmode\hbox{-}}\fi\next@}%
 \DNii@{\ifx\next-\DN@-{\leavevmode\hbox{---}}\else
  \DN@{\leavevmode\hbox{--}}\fi\next@}%
 \FN@\next@}
\def\srdr@{\kern.16667em}
\def\drsr@{\kern.02778em}
\def\sldl@{\drsr@}
\def\dlsl@{\srdr@}
\atdef@"{\unskip\relaxnext@
 \DN@{\ifx\next\space@\DN@. {\FN@\nextii@}\else
  \DN@.{\FN@\nextii@}\fi\next@.}%
 \DNii@{\ifx\next`\DN@`{\FN@\nextiii@}\else
  \ifx\next\lq\DN@\lq{\FN@\nextiii@}\else
  \DN@####1{\FN@\nextiv@}\fi\fi\next@}%
 \def\nextiii@{\ifx\next`\DN@`{\sldl@``}\else\ifx\next\lq
  \DN@\lq{\sldl@``}\else\DN@{\dlsl@`}\fi\fi\next@}%
 \def\nextiv@{\ifx\next'\DN@'{\srdr@''}\else
  \ifx\next\rq\DN@\rq{\srdr@''}\else\DN@{\drsr@'}\fi\fi\next@}%
 \FN@\next@}

\def\textfontii{\the\textfont\tw@}
\def\lbrace@{\delimiter"4266308 }
\def\rbrace@{\delimiter"5267309 }
\def\{{\RIfM@\lbrace@\else{\textfontii f}\spacefactor\@m\fi}
\def\}{\RIfM@\rbrace@\else
 \let\@sf\empty\ifhmode\edef\@sf{\spacefactor\the\spacefactor}\fi
 {\textfontii g}\@sf\relax\fi}
\let\lbrace\{
\let\rbrace\}
\def\AmSTeX{{\textfontii A\kern-.1667em%
  \lower.5ex\hbox{M}\kern-.125emS}-\TeX}
\message{line and page breaks,}
\def\vmodeerr@#1{\Err@{\string#1\space not allowed between paragraphs}}
\def\mathmodeerr@#1{\Err@{\string#1\space not allowed in math mode}}
\def\linebreak{\RIfM@\mathmodeerr@\linebreak\else
 \ifhmode\unskip\unkern\break\else\vmodeerr@\linebreak\fi\fi}

\newskip\saveskip@
\def\allowlinebreak{\RIfM@\mathmodeerr@\allowlinebreak\else
 \ifhmode\saveskip@\lastskip\unskip
 \allowbreak\ifdim\saveskip@>\z@\hskip\saveskip@\fi
 \else\vmodeerr@\allowlinebreak\fi\fi}
\def\nolinebreak{\RIfM@\mathmodeerr@\nolinebreak\else
 \ifhmode\saveskip@\lastskip\unskip
 \nobreak\ifdim\saveskip@>\z@\hskip\saveskip@\fi
 \else\vmodeerr@\nolinebreak\fi\fi}
\def\newline{\relaxnext@
 \DN@{\RIfM@\expandafter\mathmodeerr@\expandafter\newline\else
  \ifhmode\ifx\next\par\else
  \expandafter\unskip\expandafter\null\expandafter\hfill\expandafter\break\fi
  \else
  \expandafter\vmodeerr@\expandafter\newline\fi\fi}%
 \FN@\next@}
\def\dmatherr@#1{\Err@{\string#1\space not allowed in display math mode}}
\def\nondmatherr@#1{\Err@{\string#1\space not allowed in non-display math
 mode}}
\def\onlydmatherr@#1{\Err@{\string#1\space allowed only in display math mode}}
\def\nonmatherr@#1{\Err@{\string#1\space allowed only in math mode}}
\def\mathbreak{\RIfMIfI@\break\else
 \dmatherr@\mathbreak\fi\else\nonmatherr@\mathbreak\fi}
\def\nomathbreak{\RIfMIfI@\nobreak\else
 \dmatherr@\nomathbreak\fi\else\nonmatherr@\nomathbreak\fi}
\def\allowmathbreak{\RIfMIfI@\allowbreak\else
 \dmatherr@\allowmathbreak\fi\else\nonmatherr@\allowmathbreak\fi}
\def\pagebreak{\RIfM@
 \ifinner\nondmatherr@\pagebreak\else\postdisplaypenalty-\@M\fi
 \else\ifvmode\removelastskip\break\else\vadjust{\break}\fi\fi}
\def\nopagebreak{\RIfM@
 \ifinner\nondmatherr@\nopagebreak\else\postdisplaypenalty\@M\fi
 \else\ifvmode\nobreak\else\vadjust{\nobreak}\fi\fi}
\def\nonvmodeerr@#1{\Err@{\string#1\space not allowed within a paragraph
 or in math}}
\def\vnonvmode@#1#2{\relaxnext@\DNii@{\ifx\next\par\DN@{#1}\else
 \DN@{#2}\fi\next@}%
 \ifvmode\DN@{#1}\else
 \DN@{\FN@\nextii@}\fi\next@}
\def\newpage{\vnonvmode@{\vfill\break}{\nonvmodeerr@\newpage}}
\def\smallpagebreak{\vnonvmode@\smallbreak{\nonvmodeerr@\smallpagebreak}}
\def\medpagebreak{\vnonvmode@\medbreak{\nonvmodeerr@\medpagebreak}}
\def\bigpagebreak{\vnonvmode@\bigbreak{\nonvmodeerr@\bigpagebreak}}
\def\NoBlackBoxes{\global\overfullrule\z@}
\def\BlackBoxes{\global\overfullrule5\p@}
\def\Invalid@#1{\def#1{\Err@{\Invalid@@\string#1}}}
\def\Invalid@@{Invalid use of }
\message{figures,}
\Invalid@\caption
\Invalid@\captionwidth
\newdimen\smallcaptionwidth@
\def\topspace{\mid@false\ins@}
\def\midspace{\mid@true\ins@}
\newif\ifmid@
\def\captionfont@{}
\def\ins@#1{\relaxnext@\allowbreak
 \smallcaptionwidth@\captionwidth@\gdef\thespace@{#1}%
 \DN@{\ifx\next\space@\DN@. {\FN@\nextii@}\else
  \DN@.{\FN@\nextii@}\fi\next@.}%
 \DNii@{\ifx\next\caption\DN@\caption{\FN@\nextiii@}%
  \else\let\next@\nextiv@\fi\next@}%
 \def\nextiv@{\vnonvmode@
  {\ifmid@\expandafter\midinsert\else\expandafter\topinsert\fi
   \vbox to\thespace@{}\endinsert}
  {\ifmid@\nonvmodeerr@\midspace\else\nonvmodeerr@\topspace\fi}}%
 \def\nextiii@{\ifx\next\captionwidth\expandafter\nextv@
  \else\expandafter\nextvi@\fi}%
 \def\nextv@\captionwidth##1##2{\smallcaptionwidth@##1\relax\nextvi@{##2}}%
 \def\nextvi@##1{\def\thecaption@{\captionfont@##1}%
  \DN@{\ifx\next\space@\DN@. {\FN@\nextvii@}\else
   \DN@.{\FN@\nextvii@}\fi\next@.}%
  \FN@\next@}%
 \def\nextvii@{\vnonvmode@
  {\ifmid@\expandafter\midinsert\else
  \expandafter\topinsert\fi\vbox to\thespace@{}\nobreak\smallskip
  \setboxz@h{\noindent\ignorespaces\thecaption@\unskip}%
  \ifdim\wdz@>\smallcaptionwidth@\centerline{\vbox{\hsize\smallcaptionwidth@
   \noindent\ignorespaces\thecaption@\unskip}}%
  \else\centerline{\boxz@}\fi\endinsert}
  {\ifmid@\nonvmodeerr@\midspace
  \else\nonvmodeerr@\topspace\fi}}%
 \FN@\next@}
\message{comments,}
\def\newcodes@{\catcode`\\12\catcode`\{12\catcode`\}12\catcode`\#12%
 \catcode`\%12\relax}
\def\oldcodes@{\catcode`\\0\catcode`\{1\catcode`\}2\catcode`\#6%
 \catcode`\%14\relax}
\def\comment{\newcodes@\endlinechar=10 \comment@}
{\lccode`\0=`\\
\lowercase{\gdef\comment@#1^^J{\comment@@#10endcomment\comment@@@}%
\gdef\comment@@#10endcomment{\FN@\comment@@@}%
\gdef\comment@@@#1\comment@@@{\ifx\next\comment@@@\let\next\comment@
 \else\def\next{\oldcodes@\endlinechar=`\^^M\relax}%
 \fi\next}}}
\def\pr@m@s{\ifx'\next\DN@##1{\prim@s}\else\let\next@\egroup\fi\next@}
\def\prime{{\null\prime@\null}}
\mathchardef\prime@="0230
\let\dsize\displaystyle

\let\ssize\scriptstyle

\message{math spacing,}
\def\,{\RIfM@\mskip\thinmuskip\relax\else\kern.16667em\fi}
\def\!{\RIfM@\mskip-\thinmuskip\relax\else\kern-.16667em\fi}
\let\thinspace\,
\let\negthinspace\!
\def\medspace{\RIfM@\mskip\medmuskip\relax\else\kern.222222em\fi}
\def\negmedspace{\RIfM@\mskip-\medmuskip\relax\else\kern-.222222em\fi}
\def\thickspace{\RIfM@\mskip\thickmuskip\relax\else\kern.27777em\fi}
\let\;\thickspace
\def\negthickspace{\RIfM@\mskip-\thickmuskip\relax\else
 \kern-.27777em\fi}
\atdef@,{\RIfM@\mskip.1\thinmuskip\else\leavevmode\null,\fi}
\atdef@!{\RIfM@\mskip-.1\thinmuskip\else\leavevmode\null!\fi}
\atdef@.{\RIfM@&&\else\leavevmode.\spacefactor3000 \fi}
\def\and{\DOTSB\;\mathchar"3026 \;}

\message{fractions,}
\def\frac#1#2{{#1\over#2}}

\newdimen\ex@
\ex@.2326ex
\Invalid@\thickness
\def\thickfrac{\relaxnext@
 \DN@{\ifx\next\thickness\let\next@\nextii@\else
 \DN@{\nextii@\thickness1}\fi\next@}%
 \DNii@\thickness##1##2##3{{##2\above##1\ex@##3}}%
 \FN@\next@}

\def\thickfracwithdelims#1#2{\relaxnext@\def\ldelim@{#1}\def\rdelim@{#2}%
 \DN@{\ifx\next\thickness\let\next@\nextii@\else
 \DN@{\nextii@\thickness1}\fi\next@}%
 \DNii@\thickness##1##2##3{{##2\abovewithdelims
 \ldelim@\rdelim@##1\ex@##3}}%
 \FN@\next@}

\def\:{\nobreak\hskip.1111em\mathpunct{}\nonscript\mkern-\thinmuskip{:}\hskip
 .3333emplus.0555em\relax}
\def\snug{\unskip\kern-\mathsurround}
\message{smash commands,}
\def\topsmash{\top@true\bot@false\smash@}
\def\botsmash{\top@false\bot@true\smash@}
\newif\iftop@
\newif\ifbot@
\def\smash{\top@true\bot@true\smash@}
\def\smash@{\RIfM@\expandafter\mathpalette\expandafter\mathsm@sh\else
 \expandafter\makesm@sh\fi}
\def\finsm@sh{\iftop@\ht\z@\z@\fi\ifbot@\dp\z@\z@\fi\leavevmode\boxz@}
\message{large operator symbols,}
\def\LimitsOnSums{\global\let\slimits@\displaylimits}
\def\NoLimitsOnSums{\global\let\slimits@\nolimits}
\LimitsOnSums
\mathchardef\coprod@="1360       \def\coprod{\DOTSB\coprod@\slimits@}
\mathchardef\bigvee@="1357       \def\bigvee{\DOTSB\bigvee@\slimits@}
\mathchardef\bigwedge@="1356     \def\bigwedge{\DOTSB\bigwedge@\slimits@}
\mathchardef\biguplus@="1355     \def\biguplus{\DOTSB\biguplus@\slimits@}
\mathchardef\bigcap@="1354       \def\bigcap{\DOTSB\bigcap@\slimits@}
\mathchardef\bigcup@="1353       \def\bigcup{\DOTSB\bigcup@\slimits@}
\mathchardef\prod@="1351         \def\prod{\DOTSB\prod@\slimits@}
\mathchardef\sum@="1350          \def\sum{\DOTSB\sum@\slimits@}
\mathchardef\bigotimes@="134E    \def\bigotimes{\DOTSB\bigotimes@\slimits@}
\mathchardef\bigoplus@="134C     \def\bigoplus{\DOTSB\bigoplus@\slimits@}
\mathchardef\bigodot@="134A      \def\bigodot{\DOTSB\bigodot@\slimits@}
\mathchardef\bigsqcup@="1346     \def\bigsqcup{\DOTSB\bigsqcup@\slimits@}
\message{integrals,}
\def\LimitsOnInts{\global\let\ilimits@\displaylimits}
\def\NoLimitsOnInts{\global\let\ilimits@\nolimits}
\NoLimitsOnInts
\def\int{\DOTSI\intop\ilimits@}
\def\oint{\DOTSI\ointop\ilimits@}
\def\intic@{\mathchoice{\hskip.5em}{\hskip.4em}{\hskip.4em}{\hskip.4em}}
\def\negintic@{\mathchoice
 {\hskip-.5em}{\hskip-.4em}{\hskip-.4em}{\hskip-.4em}}
\def\intkern@{\mathchoice{\!\!\!}{\!\!}{\!\!}{\!\!}}
\def\intdots@{\mathchoice{\plaincdots@}
 {{\cdotp}\mkern1.5mu{\cdotp}\mkern1.5mu{\cdotp}}
 {{\cdotp}\mkern1mu{\cdotp}\mkern1mu{\cdotp}}
 {{\cdotp}\mkern1mu{\cdotp}\mkern1mu{\cdotp}}}
\newcount\intno@
\def\iint{\DOTSI\intno@\tw@\FN@\ints@}
\def\iiint{\DOTSI\intno@\thr@@\FN@\ints@}
\def\iiiint{\DOTSI\intno@4 \FN@\ints@}
\def\idotsint{\DOTSI\intno@\z@\FN@\ints@}
\def\ints@{\findlimits@\ints@@}
\newif\iflimtoken@
\newif\iflimits@
\def\findlimits@{\limtoken@true\ifx\next\limits\limits@true
 \else\ifx\next\nolimits\limits@false\else
 \limtoken@false\ifx\ilimits@\nolimits\limits@false\else
 \ifinner\limits@false\else\limits@true\fi\fi\fi\fi}
\def\multint@{\int\ifnum\intno@=\z@\intdots@                                
 \else\intkern@\fi                                                          
 \ifnum\intno@>\tw@\int\intkern@\fi                                         
 \ifnum\intno@>\thr@@\int\intkern@\fi                                       
 \int}                                                                      
\def\multintlimits@{\intop\ifnum\intno@=\z@\intdots@\else\intkern@\fi
 \ifnum\intno@>\tw@\intop\intkern@\fi
 \ifnum\intno@>\thr@@\intop\intkern@\fi\intop}
\def\ints@@{\iflimtoken@                                                    
 \def\ints@@@{\iflimits@\negintic@\mathop{\intic@\multintlimits@}\limits    
  \else\multint@\nolimits\fi                                                
  \eat@}                                                                    
 \else                                                                      
 \def\ints@@@{\iflimits@\negintic@
  \mathop{\intic@\multintlimits@}\limits\else
  \multint@\nolimits\fi}\fi\ints@@@}
\def\LimitsOnNames{\global\let\nlimits@\displaylimits}
\def\NoLimitsOnNames{\global\let\nlimits@\nolimits@}
\LimitsOnNames
\def\nolimits@{\relaxnext@
 \DN@{\ifx\next\limits\DN@\limits{\nolimits}\else
  \let\next@\nolimits\fi\next@}%
 \FN@\next@}
\message{operator names,}
\def\newmcodes@{\mathcode`\'"27\mathcode`\*"2A\mathcode`\."613A%
 \mathcode`\-"2D\mathcode`\/"2F\mathcode`\:"603A }
\def\operatorname#1{\mathop{\newmcodes@\kern\z@\fam\z@#1}\nolimits@}
\def\operatornamewithlimits#1{\mathop{\newmcodes@\kern\z@\fam\z@#1}\nlimits@}
\def\qopname@#1{\mathop{\fam\z@#1}\nolimits@}
\def\qopnamewl@#1{\mathop{\fam\z@#1}\nlimits@}
\def\arccos{\qopname@{arccos}}
\def\arcsin{\qopname@{arcsin}}
\def\arctan{\qopname@{arctan}}
\def\arg{\qopname@{arg}}
\def\cos{\qopname@{cos}}
\def\cosh{\qopname@{cosh}}
\def\cot{\qopname@{cot}}
\def\coth{\qopname@{coth}}
\def\csc{\qopname@{csc}}
\def\deg{\qopname@{deg}}
\def\det{\qopnamewl@{det}}
\def\dim{\qopname@{dim}}
\def\exp{\qopname@{exp}}
\def\gcd{\qopnamewl@{gcd}}
\def\hom{\qopname@{hom}}
\def\inf{\qopnamewl@{inf}}
\def\injlim{\qopnamewl@{inj\,lim}}
\def\ker{\qopname@{ker}}
\def\lg{\qopname@{lg}}
\def\lim{\qopnamewl@{lim}}
\def\liminf{\qopnamewl@{lim\,inf}}
\def\limsup{\qopnamewl@{lim\,sup}}
\def\ln{\qopname@{ln}}
\def\log{\qopname@{log}}
\def\max{\qopnamewl@{max}}
\def\min{\qopnamewl@{min}}
\def\Pr{\qopnamewl@{Pr}}
\def\projlim{\qopnamewl@{proj\,lim}}
\def\sec{\qopname@{sec}}
\def\sin{\qopname@{sin}}
\def\sinh{\qopname@{sinh}}
\def\sup{\qopnamewl@{sup}}
\def\tan{\qopname@{tan}}
\def\tanh{\qopname@{tanh}}
\def\varinjlim{\mathop{\vtop{\ialign{##\crcr
 \hfil\rm lim\hfil\crcr\noalign{\nointerlineskip}\rightarrowfill\crcr
 \noalign{\nointerlineskip\kern-\ex@}\crcr}}}}
\def\varprojlim{\mathop{\vtop{\ialign{##\crcr
 \hfil\rm lim\hfil\crcr\noalign{\nointerlineskip}\leftarrowfill\crcr
 \noalign{\nointerlineskip\kern-\ex@}\crcr}}}}
\def\varliminf{\mathop{\underline{\vrule height\z@ depth.2exwidth\z@
 \hbox{\rm lim}}}}

\newdimen\buffer@
\buffer@\fontdimen13 \tenex
\newdimen\buffer
\buffer\buffer@

\def\ResetBuffer{\fontdimen13 \tenex\buffer@\global\buffer\buffer@}
\def\shave#1{\mathop{\hbox{$\m@th\fontdimen13 \tenex\z@                     
 \displaystyle{#1}$}}\fontdimen13 \tenex\buffer}

\message{multilevel sub/superscripts,}
\Invalid@\\
\def\Let@{\relax\iffalse{\fi\let\\=\cr\iffalse}\fi}
\Invalid@\vspace
\def\vspace@{\def\vspace##1{\crcr\noalign{\vskip##1\relax}}}
\def\multilimits@{\bgroup\vspace@\Let@
 \baselineskip\fontdimen10 \scriptfont\tw@
 \advance\baselineskip\fontdimen12 \scriptfont\tw@
 \lineskip\thr@@\fontdimen8 \scriptfont\thr@@
 \lineskiplimit\lineskip
 \vbox\bgroup\ialign\bgroup\hfil$\m@th\scriptstyle{##}$\hfil\crcr}
\def\Sb{_\multilimits@}
\def\endSb{\crcr\egroup\egroup\egroup}
\def\Sp{^\multilimits@}

\def\spreadlines#1{\RIfMIfI@\onlydmatherr@\spreadlines\else
 \openup#1\relax\fi\else\onlydmatherr@\spreadlines\fi}
\def\Mathstrut@{\copy\Mathstrutbox@}
\newbox\Mathstrutbox@
\setbox\Mathstrutbox@\null
\setboxz@h{$\m@th($}
\ht\Mathstrutbox@\ht\z@
\dp\Mathstrutbox@\dp\z@
\message{matrices,}
\newdimen\spreadmlines@
\def\spreadmatrixlines#1{\RIfMIfI@
 \onlydmatherr@\spreadmatrixlines\else
 \spreadmlines@#1\relax\fi\else\onlydmatherr@\spreadmatrixlines\fi}
\def\matrix{\null\,\vcenter\bgroup\Let@\vspace@
 \normalbaselines\openup\spreadmlines@\ialign
 \bgroup\hfil$\m@th##$\hfil&&\quad\hfil$\m@th##$\hfil\crcr
 \Mathstrut@\crcr\noalign{\kern-\baselineskip}}
\def\endmatrix{\crcr\Mathstrut@\crcr\noalign{\kern-\baselineskip}\egroup
 \egroup\,}
\def\format{\crcr\egroup\iffalse{\fi\ifnum`}=0 \fi\format@}
\newtoks\hashtoks@
\hashtoks@{#}
\def\format@#1\\{\def\preamble@{#1}%
 \def\l{$\m@th\the\hashtoks@$\hfil}%
 \def\c{\hfil$\m@th\the\hashtoks@$\hfil}%
 \def\r{\hfil$\m@th\the\hashtoks@$}%
 \edef\preamble@@{\preamble@}\ifnum`{=0 \fi\iffalse}\fi
 \ialign\bgroup\span\preamble@@\crcr}
\def\smallmatrix{\null\,\vcenter\bgroup\vspace@\Let@
 \baselineskip9\ex@\lineskip\ex@
 \ialign\bgroup\hfil$\m@th\scriptstyle{##}$\hfil&&\thickspace\hfil
 $\m@th\scriptstyle{##}$\hfil\crcr}
\def\endsmallmatrix{\crcr\egroup\egroup\,}

\newmuskip\dotsspace@
\dotsspace@1.5mu
\def\strip@#1 {#1}
\def\spacehdots#1\for#2{\multispan{#2}\xleaders
 \hbox{$\m@th\mkern\strip@#1 \dotsspace@.\mkern\strip@#1 \dotsspace@$}\hfill}
\def\hdotsfor#1{\spacehdots\@ne\for{#1}}
\def\multispan@#1{\omit\mscount#1\unskip\loop\ifnum\mscount>\@ne\sp@n\repeat}
\def\spaceinnerhdots#1\for#2\after#3{\multispan@{\strip@#2 }#3\xleaders
 \hbox{$\m@th\mkern\strip@#1 \dotsspace@.\mkern\strip@#1 \dotsspace@$}\hfill}
\def\innerhdotsfor#1\after#2{\spaceinnerhdots\@ne\for#1\after{#2}}
\def\cases{\bgroup\spreadmlines@\jot\left\{\,\matrix\format\l&\quad\l\\}
\def\endcases{\endmatrix\right.\egroup}
\message{multiline displays,}
\newif\ifinany@
\newif\ifinalign@
\newif\ifingather@
\def\strut@{\copy\strutbox@}
\newbox\strutbox@
\setbox\strutbox@\hbox{\vrule height8\p@ depth3\p@ width\z@}
\def\topaligned{\null\,\vtop\aligned@}
\def\botaligned{\null\,\vbox\aligned@}
\def\aligned{\null\,\vcenter\aligned@}
\def\aligned@{\bgroup\vspace@\Let@
 \ifinany@\else\openup\jot\fi\ialign
 \bgroup\hfil\strut@$\m@th\displaystyle{##}$&
 $\m@th\displaystyle{{}##}$\hfil\crcr}
\def\endaligned{\crcr\egroup\egroup}

\def\alignedat#1{\null\,\vcenter\bgroup\doat@{#1}\vspace@\Let@
 \ifinany@\else\openup\jot\fi\ialign\bgroup\span\preamble@@\crcr}
\newcount\atcount@
\def\doat@#1{\toks@{\hfil\strut@$\m@th
 \displaystyle{\the\hashtoks@}$&$\m@th\displaystyle
 {{}\the\hashtoks@}$\hfil}
 \atcount@#1\relax\advance\atcount@\m@ne                                    
 \loop\ifnum\atcount@>\z@\toks@=\expandafter{\the\toks@&\hfil$\m@th
 \displaystyle{\the\hashtoks@}$&$\m@th
 \displaystyle{{}\the\hashtoks@}$\hfil}\advance
  \atcount@\m@ne\repeat                                                     
 \xdef\preamble@{\the\toks@}\xdef\preamble@@{\preamble@}}

\def\gathered{\null\,\vcenter\bgroup\vspace@\Let@
 \ifinany@\else\openup\jot\fi\ialign
 \bgroup\hfil\strut@$\m@th\displaystyle{##}$\hfil\crcr}
\def\endgathered{\crcr\egroup\egroup}
\newif\iftagsleft@
\def\TagsOnLeft{\global\tagsleft@true}
\def\TagsOnRight{\global\tagsleft@false}
\TagsOnLeft
\newif\ifmathtags@
\def\TagsAsMath{\global\mathtags@true}
\def\TagsAsText{\global\mathtags@false}
\TagsAsText
\def\tagform@#1{\hbox{\rm(\ignorespaces#1\unskip)}}
\def\thetag{\leavevmode\tagform@}
\def\tag#1$${\iftagsleft@\leqno\else\eqno\fi                                
 \maketag@#1\maketag@                                                       
 $$}                                                                        
\def\maketag@{\FN@\maketag@@}
\def\maketag@@{\ifx\next"\expandafter\maketag@@@\else\expandafter\maketag@@@@
 \fi}
\def\maketag@@@"#1"#2\maketag@{\hbox{\rm#1}}                                
\def\maketag@@@@#1\maketag@{\ifmathtags@\tagform@{$\m@th#1$}\else
 \tagform@{#1}\fi}
\interdisplaylinepenalty\@M
\def\allowdisplaybreaks{\RIfMIfI@
 \onlydmatherr@\allowdisplaybreaks\else
 \interdisplaylinepenalty\z@\fi\else\onlydmatherr@\allowdisplaybreaks\fi}
\Invalid@\allowdisplaybreak
\Invalid@\displaybreak
\Invalid@\intertext
\def\allowdisplaybreak@{\def\allowdisplaybreak{\crcr\noalign{\allowbreak}}}
\def\displaybreak@{\def\displaybreak{\crcr\noalign{\break}}}
\def\intertext@{\def\intertext##1{\crcr\noalign{%
 \penalty\postdisplaypenalty \vskip\belowdisplayskip
 \vbox{\normalbaselines\noindent##1}%
 \penalty\predisplaypenalty \vskip\abovedisplayskip}}}
\newskip\centering@
\centering@\z@ plus\@m\p@
\def\align{\relax\ifingather@\DN@{\csname align (in
  \string\gather)\endcsname}\else
 \ifmmode\ifinner\DN@{\onlydmatherr@\align}\else
  \let\next@\align@\fi
 \else\DN@{\onlydmatherr@\align}\fi\fi\next@}
\newhelp\andhelp@
{An extra & here is so disastrous that you should probably exit^^J
and fix things up.}
\newif\iftag@
\newcount\and@
\def\align@{\inalign@true\inany@true
 \vspace@\allowdisplaybreak@\displaybreak@\intertext@
 \def\tag{\global\tag@true\ifnum\and@=\z@\DN@{&&}\else
          \DN@{&}\fi\next@}%
 \iftagsleft@\DN@{\csname align \endcsname}\else
  \DN@{\csname align \space\endcsname}\fi\next@}
\def\Tag@{\iftag@\else\errhelp\andhelp@\err@{Extra & on this line}\fi}
\newdimen\lwidth@
\newdimen\rwidth@
\newdimen\maxlwidth@
\newdimen\maxrwidth@
\newdimen\totwidth@
\def\measure@#1\endalign{\lwidth@\z@\rwidth@\z@\maxlwidth@\z@\maxrwidth@\z@
 \global\and@\z@                                                            
 \setbox@ne\vbox                                                            
  {\everycr{\noalign{\global\tag@false\global\and@\z@}}\Let@                
  \halign{\setboxz@h{$\m@th\displaystyle{\@lign##}$}
   \global\lwidth@\wdz@                                                     
   \ifdim\lwidth@>\maxlwidth@\global\maxlwidth@\lwidth@\fi                  
   \global\advance\and@\@ne                                                 
   &\setboxz@h{$\m@th\displaystyle{{}\@lign##}$}\global\rwidth@\wdz@        
   \ifdim\rwidth@>\maxrwidth@\global\maxrwidth@\rwidth@\fi                  
   \global\advance\and@\@ne                                                
   &\Tag@
   \eat@{##}\crcr#1\crcr}}
 \totwidth@\maxlwidth@\advance\totwidth@\maxrwidth@}                       
\def\displ@y@{\global\dt@ptrue\openup\jot
 \everycr{\noalign{\global\tag@false\global\and@\z@\ifdt@p\global\dt@pfalse
 \vskip-\lineskiplimit\vskip\normallineskiplimit\else
 \penalty\interdisplaylinepenalty\fi}}}
\def\black@#1{\noalign{\ifdim#1>\displaywidth
 \dimen@\prevdepth\nointerlineskip                                          
 \vskip-\ht\strutbox@\vskip-\dp\strutbox@                                   
 \vbox{\noindent\hbox to#1{\strut@\hfill}}
 \prevdepth\dimen@                                                          
 \fi}}
\expandafter\def\csname align \space\endcsname#1\endalign
 {\measure@#1\endalign\global\and@\z@                                       
 \ifingather@\everycr{\noalign{\global\and@\z@}}\else\displ@y@\fi           
 \Let@\tabskip\centering@                                                   
 \halign to\displaywidth
  {\hfil\strut@\setboxz@h{$\m@th\displaystyle{\@lign##}$}
  \global\lwidth@\wdz@\boxz@\global\advance\and@\@ne                        
  \tabskip\z@skip                                                           
  &\setboxz@h{$\m@th\displaystyle{{}\@lign##}$}
  \global\rwidth@\wdz@\boxz@\hfill\global\advance\and@\@ne                  
  \tabskip\centering@                                                       
  &\setboxz@h{\@lign\strut@\maketag@##\maketag@}
  \dimen@\displaywidth\advance\dimen@-\totwidth@
  \divide\dimen@\tw@\advance\dimen@\maxrwidth@\advance\dimen@-\rwidth@     
  \ifdim\dimen@<\tw@\wdz@\llap{\vtop{\normalbaselines\null\boxz@}}
  \else\llap{\boxz@}\fi                                                    
  \tabskip\z@skip                                                          
  \crcr#1\crcr                                                             
  \black@\totwidth@}}                                                      
\newdimen\lineht@
\expandafter\def\csname align \endcsname#1\endalign{\measure@#1\endalign
 \global\and@\z@
 \ifdim\totwidth@>\displaywidth\let\displaywidth@\totwidth@\else
  \let\displaywidth@\displaywidth\fi                                        
 \ifingather@\everycr{\noalign{\global\and@\z@}}\else\displ@y@\fi
 \Let@\tabskip\centering@\halign to\displaywidth
  {\hfil\strut@\setboxz@h{$\m@th\displaystyle{\@lign##}$}%
  \global\lwidth@\wdz@\global\lineht@\ht\z@                                 
  \boxz@\global\advance\and@\@ne
  \tabskip\z@skip&\setboxz@h{$\m@th\displaystyle{{}\@lign##}$}%
  \global\rwidth@\wdz@\ifdim\ht\z@>\lineht@\global\lineht@\ht\z@\fi         
  \boxz@\hfil\global\advance\and@\@ne
  \tabskip\centering@&\kern-\displaywidth@                                  
  \setboxz@h{\@lign\strut@\maketag@##\maketag@}%
  \dimen@\displaywidth\advance\dimen@-\totwidth@
  \divide\dimen@\tw@\advance\dimen@\maxlwidth@\advance\dimen@-\lwidth@
  \ifdim\dimen@<\tw@\wdz@
   \rlap{\vbox{\normalbaselines\boxz@\vbox to\lineht@{}}}\else
   \rlap{\boxz@}\fi
  \tabskip\displaywidth@\crcr#1\crcr\black@\totwidth@}}
\expandafter\def\csname align (in \string\gather)\endcsname
  #1\endalign{\vcenter{\align@#1\endalign}}
\Invalid@\endalign
\newif\ifxat@
\def\alignat{\RIfMIfI@\DN@{\onlydmatherr@\alignat}\else
 \DN@{\csname alignat \endcsname}\fi\else
 \DN@{\onlydmatherr@\alignat}\fi\next@}
\newif\ifmeasuring@
\newbox\savealignat@
\expandafter\def\csname alignat \endcsname#1#2\endalignat                   
 {\inany@true\xat@false
 \def\tag{\global\tag@true\count@#1\relax\multiply\count@\tw@
  \xdef\tag@{}\loop\ifnum\count@>\and@\xdef\tag@{&\tag@}\advance\count@\m@ne
  \repeat\tag@}%
 \vspace@\allowdisplaybreak@\displaybreak@\intertext@
 \displ@y@\measuring@true                                                   
 \setbox\savealignat@\hbox{$\m@th\displaystyle\Let@
  \attag@{#1}
  \vbox{\halign{\span\preamble@@\crcr#2\crcr}}$}%
 \measuring@false                                                           
 \Let@\attag@{#1}
 \tabskip\centering@\halign to\displaywidth
  {\span\preamble@@\crcr#2\crcr                                             
  \black@{\wd\savealignat@}}}                                               
\Invalid@\endalignat
\def\xalignat{\RIfMIfI@
 \DN@{\onlydmatherr@\xalignat}\else
 \DN@{\csname xalignat \endcsname}\fi\else
 \DN@{\onlydmatherr@\xalignat}\fi\next@}
\expandafter\def\csname xalignat \endcsname#1#2\endxalignat
 {\inany@true\xat@true
 \def\tag{\global\tag@true\def\tag@{}\count@#1\relax\multiply\count@\tw@
  \loop\ifnum\count@>\and@\xdef\tag@{&\tag@}\advance\count@\m@ne\repeat\tag@}%
 \vspace@\allowdisplaybreak@\displaybreak@\intertext@
 \displ@y@\measuring@true\setbox\savealignat@\hbox{$\m@th\displaystyle\Let@
 \attag@{#1}\vbox{\halign{\span\preamble@@\crcr#2\crcr}}$}%
 \measuring@false\Let@
 \attag@{#1}\tabskip\centering@\halign to\displaywidth
 {\span\preamble@@\crcr#2\crcr\black@{\wd\savealignat@}}}
\def\attag@#1{\let\Maketag@\maketag@\let\TAG@\Tag@                          
 \let\Tag@=0\let\maketag@=0
 \ifmeasuring@\def\llap@##1{\setboxz@h{##1}\hbox to\tw@\wdz@{}}%
  \def\rlap@##1{\setboxz@h{##1}\hbox to\tw@\wdz@{}}\else
  \let\llap@\llap\let\rlap@\rlap\fi                                         
 \toks@{\hfil\strut@$\m@th\displaystyle{\@lign\the\hashtoks@}$\tabskip\z@skip
  \global\advance\and@\@ne&$\m@th\displaystyle{{}\@lign\the\hashtoks@}$\hfil
  \ifxat@\tabskip\centering@\fi\global\advance\and@\@ne}
 \iftagsleft@
  \toks@@{\tabskip\centering@&\Tag@\kern-\displaywidth
   \rlap@{\@lign\maketag@\the\hashtoks@\maketag@}%
   \global\advance\and@\@ne\tabskip\displaywidth}\else
  \toks@@{\tabskip\centering@&\Tag@\llap@{\@lign\maketag@
   \the\hashtoks@\maketag@}\global\advance\and@\@ne\tabskip\z@skip}\fi      
 \atcount@#1\relax\advance\atcount@\m@ne
 \loop\ifnum\atcount@>\z@
 \toks@=\expandafter{\the\toks@&\hfil$\m@th\displaystyle{\@lign
  \the\hashtoks@}$\global\advance\and@\@ne
  \tabskip\z@skip&$\m@th\displaystyle{{}\@lign\the\hashtoks@}$\hfil\ifxat@
  \tabskip\centering@\fi\global\advance\and@\@ne}\advance\atcount@\m@ne
 \repeat                                                                    
 \xdef\preamble@{\the\toks@\the\toks@@}
 \xdef\preamble@@{\preamble@}
 \let\maketag@\Maketag@\let\Tag@\TAG@}                                      
\Invalid@\endxalignat
\def\xxalignat{\RIfMIfI@
 \DN@{\onlydmatherr@\xxalignat}\else\DN@{\csname xxalignat
  \endcsname}\fi\else
 \DN@{\onlydmatherr@\xxalignat}\fi\next@}
\expandafter\def\csname xxalignat \endcsname#1#2\endxxalignat{\inany@true
 \vspace@\allowdisplaybreak@\displaybreak@\intertext@
 \displ@y\setbox\savealignat@\hbox{$\m@th\displaystyle\Let@
 \xxattag@{#1}\vbox{\halign{\span\preamble@@\crcr#2\crcr}}$}%
 \Let@\xxattag@{#1}\tabskip\z@skip\halign to\displaywidth
 {\span\preamble@@\crcr#2\crcr\black@{\wd\savealignat@}}}
\def\xxattag@#1{\toks@{\tabskip\z@skip\hfil\strut@
 $\m@th\displaystyle{\the\hashtoks@}$&%
 $\m@th\displaystyle{{}\the\hashtoks@}$\hfil\tabskip\centering@&}%
 \atcount@#1\relax\advance\atcount@\m@ne\loop\ifnum\atcount@>\z@
 \toks@=\expandafter{\the\toks@&\hfil$\m@th\displaystyle{\the\hashtoks@}$%
  \tabskip\z@skip&$\m@th\displaystyle{{}\the\hashtoks@}$\hfil
  \tabskip\centering@}\advance\atcount@\m@ne\repeat
 \xdef\preamble@{\the\toks@\tabskip\z@skip}\xdef\preamble@@{\preamble@}}
\Invalid@\endxxalignat
\newdimen\gwidth@
\newdimen\gmaxwidth@
\def\gmeasure@#1\endgather{\gwidth@\z@\gmaxwidth@\z@\setbox@ne\vbox{\Let@
 \halign{\setboxz@h{$\m@th\displaystyle{##}$}\global\gwidth@\wdz@
 \ifdim\gwidth@>\gmaxwidth@\global\gmaxwidth@\gwidth@\fi
 &\eat@{##}\crcr#1\crcr}}}
\def\gather{\RIfMIfI@\DN@{\onlydmatherr@\gather}\else
 \ingather@true\inany@true\def\tag{&}%
 \vspace@\allowdisplaybreak@\displaybreak@\intertext@
 \displ@y\Let@
 \iftagsleft@\DN@{\csname gather \endcsname}\else
  \DN@{\csname gather \space\endcsname}\fi\fi
 \else\DN@{\onlydmatherr@\gather}\fi\next@}
\expandafter\def\csname gather \space\endcsname#1\endgather
 {\gmeasure@#1\endgather\tabskip\centering@
 \halign to\displaywidth{\hfil\strut@\setboxz@h{$\m@th\displaystyle{##}$}%
 \global\gwidth@\wdz@\boxz@\hfil&
 \setboxz@h{\strut@{\maketag@##\maketag@}}%
 \dimen@\displaywidth\advance\dimen@-\gwidth@
 \ifdim\dimen@>\tw@\wdz@\llap{\boxz@}\else
 \llap{\vtop{\normalbaselines\null\boxz@}}\fi
 \tabskip\z@skip\crcr#1\crcr\black@\gmaxwidth@}}
\newdimen\glineht@
\expandafter\def\csname gather \endcsname#1\endgather{\gmeasure@#1\endgather
 \ifdim\gmaxwidth@>\displaywidth\let\gdisplaywidth@\gmaxwidth@\else
 \let\gdisplaywidth@\displaywidth\fi\tabskip\centering@\halign to\displaywidth
 {\hfil\strut@\setboxz@h{$\m@th\displaystyle{##}$}%
 \global\gwidth@\wdz@\global\glineht@\ht\z@\boxz@\hfil&\kern-\gdisplaywidth@
 \setboxz@h{\strut@{\maketag@##\maketag@}}%
 \dimen@\displaywidth\advance\dimen@-\gwidth@
 \ifdim\dimen@>\tw@\wdz@\rlap{\boxz@}\else
 \rlap{\vbox{\normalbaselines\boxz@\vbox to\glineht@{}}}\fi
 \tabskip\gdisplaywidth@\crcr#1\crcr\black@\gmaxwidth@}}
\newif\ifctagsplit@
\def\CenteredTagsOnSplits{\global\ctagsplit@true}
\def\TopOrBottomTagsOnSplits{\global\ctagsplit@false}
\TopOrBottomTagsOnSplits
\def\split{\relax\ifinany@\let\next@\insplit@\else
 \ifmmode\ifinner\def\next@{\onlydmatherr@\split}\else
 \let\next@\outsplit@\fi\else
 \def\next@{\onlydmatherr@\split}\fi\fi\next@}
\def\insplit@{\global\setbox\z@\vbox\bgroup\vspace@\Let@\ialign\bgroup
 \hfil\strut@$\m@th\displaystyle{##}$&$\m@th\displaystyle{{}##}$\hfill\crcr}
\def\endsplit{\crcr\egroup\egroup\iftagsleft@\expandafter\lendsplit@\else
 \expandafter\rendsplit@\fi}
\def\rendsplit@{\global\setbox9 \vbox
 {\unvcopy\z@\global\setbox8 \lastbox\unskip}
 \setbox@ne\hbox{\unhcopy8 \unskip\global\setbox\tw@\lastbox
 \unskip\global\setbox\thr@@\lastbox}
 \global\setbox7 \hbox{\unhbox\tw@\unskip}
 \ifinalign@\ifctagsplit@                                                   
  \gdef\split@{\hbox to\wd\thr@@{}&
   \vcenter{\vbox{\moveleft\wd\thr@@\boxz@}}}
 \else\gdef\split@{&\vbox{\moveleft\wd\thr@@\box9}\crcr
  \box\thr@@&\box7}\fi                                                      
 \else                                                                      
  \ifctagsplit@\gdef\split@{\vcenter{\boxz@}}\else
  \gdef\split@{\box9\crcr\hbox{\box\thr@@\box7}}\fi
 \fi
 \split@}                                                                   
\def\lendsplit@{\global\setbox9\vtop{\unvcopy\z@}
 \setbox@ne\vbox{\unvcopy\z@\global\setbox8\lastbox}
 \setbox@ne\hbox{\unhcopy8\unskip\setbox\tw@\lastbox
  \unskip\global\setbox\thr@@\lastbox}
 \ifinalign@\ifctagsplit@                                                   
  \gdef\split@{\hbox to\wd\thr@@{}&
  \vcenter{\vbox{\moveleft\wd\thr@@\box9}}}
  \else                                                                     
  \gdef\split@{\hbox to\wd\thr@@{}&\vbox{\moveleft\wd\thr@@\box9}}\fi
 \else
  \ifctagsplit@\gdef\split@{\vcenter{\box9}}\else
  \gdef\split@{\box9}\fi
 \fi\split@}
\def\outsplit@#1$${\align\insplit@#1\endalign$$}
\newdimen\multlinegap@
\multlinegap@1em
\newdimen\multlinetaggap@
\multlinetaggap@1em
\def\MultlineGap#1{\global\multlinegap@#1\relax}
\def\multlinegap#1{\RIfMIfI@\onlydmatherr@\multlinegap\else
 \multlinegap@#1\relax\fi\else\onlydmatherr@\multlinegap\fi}
\def\nomultlinegap{\multlinegap{\z@}}
\def\multline{\RIfMIfI@
 \DN@{\onlydmatherr@\multline}\else
 \DN@{\multline@}\fi\else
 \DN@{\onlydmatherr@\multline}\fi\next@}
\newif\iftagin@
\def\tagin@#1{\tagin@false\in@\tag{#1}\ifin@\tagin@true\fi}
\def\multline@#1$${\inany@true\vspace@\allowdisplaybreak@\displaybreak@
 \tagin@{#1}\iftagsleft@\DN@{\multline@l#1$$}\else
 \DN@{\multline@r#1$$}\fi\next@}
\newdimen\mwidth@
\def\rmmeasure@#1\endmultline{%
 \def\shoveleft##1{##1}\def\shoveright##1{##1}
 \setbox@ne\vbox{\Let@\halign{\setboxz@h
  {$\m@th\@lign\displaystyle{}##$}\global\mwidth@\wdz@
  \crcr#1\crcr}}}
\newdimen\mlineht@
\newif\ifzerocr@
\newif\ifonecr@
\def\lmmeasure@#1\endmultline{\global\zerocr@true\global\onecr@false
 \everycr{\noalign{\ifonecr@\global\onecr@false\fi
  \ifzerocr@\global\zerocr@false\global\onecr@true\fi}}
  \def\shoveleft##1{##1}\def\shoveright##1{##1}%
 \setbox@ne\vbox{\Let@\halign{\setboxz@h
  {$\m@th\@lign\displaystyle{}##$}\ifonecr@\global\mwidth@\wdz@
  \global\mlineht@\ht\z@\fi\crcr#1\crcr}}}
\newbox\mtagbox@
\newdimen\ltwidth@
\newdimen\rtwidth@
\def\multline@l#1$${\iftagin@\DN@{\lmultline@@#1$$}\else
 \DN@{\setbox\mtagbox@\null\ltwidth@\z@\rtwidth@\z@
  \lmultline@@@#1$$}\fi\next@}
\def\lmultline@@#1\endmultline\tag#2$${%
 \setbox\mtagbox@\hbox{\maketag@#2\maketag@}
 \lmmeasure@#1\endmultline\dimen@\mwidth@\advance\dimen@\wd\mtagbox@
 \advance\dimen@\multlinetaggap@                                            
 \ifdim\dimen@>\displaywidth\ltwidth@\z@\else\ltwidth@\wd\mtagbox@\fi       
 \lmultline@@@#1\endmultline$$}
\def\lmultline@@@{\displ@y
 \def\shoveright##1{##1\hfilneg\hskip\multlinegap@}%
 \def\shoveleft##1{\setboxz@h{$\m@th\displaystyle{}##1$}%
  \setbox@ne\hbox{$\m@th\displaystyle##1$}%
  \hfilneg
  \iftagin@
   \ifdim\ltwidth@>\z@\hskip\ltwidth@\hskip\multlinetaggap@\fi
  \else\hskip\multlinegap@\fi\hskip.5\wd@ne\hskip-.5\wdz@##1}
  \halign\bgroup\Let@\hbox to\displaywidth
   {\strut@$\m@th\displaystyle\hfil{}##\hfil$}\crcr
   \hfilneg                                                                 
   \iftagin@                                                                
    \ifdim\ltwidth@>\z@                                                     
     \box\mtagbox@\hskip\multlinetaggap@                                    
    \else
     \rlap{\vbox{\normalbaselines\hbox{\strut@\box\mtagbox@}%
     \vbox to\mlineht@{}}}\fi                                               
   \else\hskip\multlinegap@\fi}                                             
\def\multline@r#1$${\iftagin@\DN@{\rmultline@@#1$$}\else
 \DN@{\setbox\mtagbox@\null\ltwidth@\z@\rtwidth@\z@
  \rmultline@@@#1$$}\fi\next@}
\def\rmultline@@#1\endmultline\tag#2$${\ltwidth@\z@
 \setbox\mtagbox@\hbox{\maketag@#2\maketag@}%
 \rmmeasure@#1\endmultline\dimen@\mwidth@\advance\dimen@\wd\mtagbox@
 \advance\dimen@\multlinetaggap@
 \ifdim\dimen@>\displaywidth\rtwidth@\z@\else\rtwidth@\wd\mtagbox@\fi
 \rmultline@@@#1\endmultline$$}
\def\rmultline@@@{\displ@y
 \def\shoveright##1{##1\hfilneg\iftagin@\ifdim\rtwidth@>\z@
  \hskip\rtwidth@\hskip\multlinetaggap@\fi\else\hskip\multlinegap@\fi}%
 \def\shoveleft##1{\setboxz@h{$\m@th\displaystyle{}##1$}%
  \setbox@ne\hbox{$\m@th\displaystyle##1$}%
  \hfilneg\hskip\multlinegap@\hskip.5\wd@ne\hskip-.5\wdz@##1}%
 \halign\bgroup\Let@\hbox to\displaywidth
  {\strut@$\m@th\displaystyle\hfil{}##\hfil$}\crcr
 \hfilneg\hskip\multlinegap@}
\def\endmultline{\iftagsleft@\expandafter\lendmultline@\else
 \expandafter\rendmultline@\fi}
\def\lendmultline@{\hfilneg\hskip\multlinegap@\crcr\egroup}
\def\rendmultline@{\iftagin@                                                
 \ifdim\rtwidth@>\z@                                                        
  \hskip\multlinetaggap@\box\mtagbox@                                       
 \else\llap{\vtop{\normalbaselines\null\hbox{\strut@\box\mtagbox@}}}\fi     
 \else\hskip\multlinegap@\fi                                                
 \hfilneg\crcr\egroup}
\def\bmod{\mskip-\medmuskip\mkern5mu\mathbin{\fam\z@ mod}\penalty900
 \mkern5mu\mskip-\medmuskip}
\def\pmod#1{\allowbreak\ifinner\mkern8mu\else\mkern18mu\fi
 ({\fam\z@ mod}\,\,#1)}
\def\pod#1{\allowbreak\ifinner\mkern8mu\else\mkern18mu\fi(#1)}
\def\mod#1{\allowbreak\ifinner\mkern12mu\else\mkern18mu\fi{\fam\z@ mod}\,\,#1}
\message{continued fractions,}
\newcount\cfraccount@
\def\cfrac{\bgroup\bgroup\advance\cfraccount@\@ne\strut
 \iffalse{\fi\def\\{\over\displaystyle}\iffalse}\fi}
\def\lcfrac{\bgroup\bgroup\advance\cfraccount@\@ne\strut
 \iffalse{\fi\def\\{\hfill\over\displaystyle}\iffalse}\fi}
\def\rcfrac{\bgroup\bgroup\advance\cfraccount@\@ne\strut\hfill
 \iffalse{\fi\def\\{\over\displaystyle}\iffalse}\fi}
\def\gloop@#1\repeat{\gdef\body{#1}\iterate}
\def\endcfrac{\gloop@\ifnum\cfraccount@>\z@\global\advance\cfraccount@\m@ne
 \egroup\hskip-\nulldelimiterspace\egroup\repeat}
\message{compound symbols,}
\def\binrel@#1{\setboxz@h{\thinmuskip0mu
  \medmuskip\m@ne mu\thickmuskip\@ne mu$#1\m@th$}%
 \setbox@ne\hbox{\thinmuskip0mu\medmuskip\m@ne mu\thickmuskip
  \@ne mu${}#1{}\m@th$}%
 \setbox\tw@\hbox{\hskip\wd@ne\hskip-\wdz@}}
\def\overset#1\to#2{\binrel@{#2}\ifdim\wd\tw@<\z@
 \mathbin{\mathop{\kern\z@#2}\limits^{#1}}\else\ifdim\wd\tw@>\z@
 \mathrel{\mathop{\kern\z@#2}\limits^{#1}}\else
 {\mathop{\kern\z@#2}\limits^{#1}}{}\fi\fi}
\def\underset#1\to#2{\binrel@{#2}\ifdim\wd\tw@<\z@
 \mathbin{\mathop{\kern\z@#2}\limits_{#1}}\else\ifdim\wd\tw@>\z@
 \mathrel{\mathop{\kern\z@#2}\limits_{#1}}\else
 {\mathop{\kern\z@#2}\limits_{#1}}{}\fi\fi}
\def\oversetbrace#1\to#2{\overbrace{#2}^{#1}}
\def\undersetbrace#1\to#2{\underbrace{#2}_{#1}}
\def\sideset#1\and#2\to#3{%
 \setbox@ne\hbox{$\dsize{\vphantom{#3}}#1{#3}\m@th$}%
 \setbox\tw@\hbox{$\dsize{#3}#2\m@th$}%
 \hskip\wd@ne\hskip-\wd\tw@\mathop{\hskip\wd\tw@\hskip-\wd@ne
  {\vphantom{#3}}#1{#3}#2}}
\def\rightarrowfill@#1{\setboxz@h{$#1-\m@th$}\ht\z@\z@
  $#1\m@th\copy\z@\mkern-6mu\cleaders
  \hbox{$#1\mkern-2mu\box\z@\mkern-2mu$}\hfill
  \mkern-6mu\mathord\rightarrow$}
\def\leftarrowfill@#1{\setboxz@h{$#1-\m@th$}\ht\z@\z@
  $#1\m@th\mathord\leftarrow\mkern-6mu\cleaders
  \hbox{$#1\mkern-2mu\copy\z@\mkern-2mu$}\hfill
  \mkern-6mu\box\z@$}
\def\leftrightarrowfill@#1{\setboxz@h{$#1-\m@th$}\ht\z@\z@
  $#1\m@th\mathord\leftarrow\mkern-6mu\cleaders
  \hbox{$#1\mkern-2mu\box\z@\mkern-2mu$}\hfill
  \mkern-6mu\mathord\rightarrow$}
\def\overrightarrow{\mathpalette\overrightarrow@}
\def\overrightarrow@#1#2{\vbox{\ialign{##\crcr\rightarrowfill@#1\crcr
 \noalign{\kern-\ex@\nointerlineskip}$\m@th\hfil#1#2\hfil$\crcr}}}

\def\overleftarrow{\mathpalette\overleftarrow@}
\def\overleftarrow@#1#2{\vbox{\ialign{##\crcr\leftarrowfill@#1\crcr
 \noalign{\kern-\ex@\nointerlineskip}$\m@th\hfil#1#2\hfil$\crcr}}}
\def\overleftrightarrow{\mathpalette\overleftrightarrow@}
\def\overleftrightarrow@#1#2{\vbox{\ialign{##\crcr\leftrightarrowfill@#1\crcr
 \noalign{\kern-\ex@\nointerlineskip}$\m@th\hfil#1#2\hfil$\crcr}}}
\def\underrightarrow{\mathpalette\underrightarrow@}
\def\underrightarrow@#1#2{\vtop{\ialign{##\crcr$\m@th\hfil#1#2\hfil$\crcr
 \noalign{\nointerlineskip}\rightarrowfill@#1\crcr}}}

\def\underleftarrow{\mathpalette\underleftarrow@}
\def\underleftarrow@#1#2{\vtop{\ialign{##\crcr$\m@th\hfil#1#2\hfil$\crcr
 \noalign{\nointerlineskip}\leftarrowfill@#1\crcr}}}
\def\underleftrightarrow{\mathpalette\underleftrightarrow@}
\def\underleftrightarrow@#1#2{\vtop{\ialign{##\crcr$\m@th\hfil#1#2\hfil$\crcr
 \noalign{\nointerlineskip}\leftrightarrowfill@#1\crcr}}}
\message{various kinds of dots,}
\let\DOTSI\relax
\let\DOTSB\relax

\newif\ifmath@
{\uccode`7=`\\ \uccode`8=`m \uccode`9=`a \uccode`0=`t \uccode`!=`h
 \uppercase{\gdef\math@#1#2#3#4#5#6\math@{\global\math@false\ifx 7#1\ifx 8#2%
 \ifx 9#3\ifx 0#4\ifx !#5\xdef\meaning@{#6}\global\math@true\fi\fi\fi\fi\fi}}}
\newif\ifmathch@
{\uccode`7=`c \uccode`8=`h \uccode`9=`\"
 \uppercase{\gdef\mathch@#1#2#3#4#5#6\mathch@{\global\mathch@false
  \ifx 7#1\ifx 8#2\ifx 9#5\global\mathch@true\xdef\meaning@{9#6}\fi\fi\fi}}}
\newcount\classnum@
\def\getmathch@#1.#2\getmathch@{\classnum@#1 \divide\classnum@4096
 \ifcase\number\classnum@\or\or\gdef\thedots@{\dotsb@}\or
 \gdef\thedots@{\dotsb@}\fi}
\newif\ifmathbin@
{\uccode`4=`b \uccode`5=`i \uccode`6=`n
 \uppercase{\gdef\mathbin@#1#2#3{\relaxnext@
  \DNii@##1\mathbin@{\ifx\space@\next\global\mathbin@true\fi}%
 \global\mathbin@false\DN@##1\mathbin@{}%
 \ifx 4#1\ifx 5#2\ifx 6#3\DN@{\FN@\nextii@}\fi\fi\fi\next@}}}
\newif\ifmathrel@
{\uccode`4=`r \uccode`5=`e \uccode`6=`l
 \uppercase{\gdef\mathrel@#1#2#3{\relaxnext@
  \DNii@##1\mathrel@{\ifx\space@\next\global\mathrel@true\fi}%
 \global\mathrel@false\DN@##1\mathrel@{}%
 \ifx 4#1\ifx 5#2\ifx 6#3\DN@{\FN@\nextii@}\fi\fi\fi\next@}}}
\newif\ifmacro@
{\uccode`5=`m \uccode`6=`a \uccode`7=`c
 \uppercase{\gdef\macro@#1#2#3#4\macro@{\global\macro@false
  \ifx 5#1\ifx 6#2\ifx 7#3\global\macro@true
  \xdef\meaning@{\macro@@#4\macro@@}\fi\fi\fi}}}
\def\macro@@#1->#2\macro@@{#2}
\newif\ifDOTS@
\newcount\DOTSCASE@
{\uccode`6=`\\ \uccode`7=`D \uccode`8=`O \uccode`9=`T \uccode`0=`S
 \uppercase{\gdef\DOTS@#1#2#3#4#5{\global\DOTS@false\DN@##1\DOTS@{}%
  \ifx 6#1\ifx 7#2\ifx 8#3\ifx 9#4\ifx 0#5\let\next@\DOTS@@\fi\fi\fi\fi\fi
  \next@}}}
{\uccode`3=`B \uccode`4=`I \uccode`5=`X
 \uppercase{\gdef\DOTS@@#1{\relaxnext@
  \DNii@##1\DOTS@{\ifx\space@\next\global\DOTS@true\fi}%
  \DN@{\FN@\nextii@}%
  \ifx 3#1\global\DOTSCASE@\z@\else
  \ifx 4#1\global\DOTSCASE@\@ne\else
  \ifx 5#1\global\DOTSCASE@\tw@\else\DN@##1\DOTS@{}%
  \fi\fi\fi\next@}}}
\newif\ifnot@
{\uccode`5=`\\ \uccode`6=`n \uccode`7=`o \uccode`8=`t
 \uppercase{\gdef\not@#1#2#3#4{\relaxnext@
  \DNii@##1\not@{\ifx\space@\next\global\not@true\fi}%
 \global\not@false\DN@##1\not@{}%
 \ifx 5#1\ifx 6#2\ifx 7#3\ifx 8#4\DN@{\FN@\nextii@}\fi\fi\fi
 \fi\next@}}}
\newif\ifkeybin@
\def\keybin@{\keybin@true
 \ifx\next+\else\ifx\next=\else\ifx\next<\else\ifx\next>\else\ifx\next-\else
 \ifx\next*\else\ifx\next:\else\keybin@false\fi\fi\fi\fi\fi\fi\fi}
\def\dots{\RIfM@\expandafter\mdots@\else\expandafter\tdots@\fi}
\def\tdots@{\unskip\relaxnext@
 \DN@{$\m@th\mathinner{\ldotp\ldotp\ldotp}\,
   \ifx\next,\,$\else\ifx\next.\,$\else\ifx\next;\,$\else\ifx\next:\,$\else
   \ifx\next?\,$\else\ifx\next!\,$\else$ \fi\fi\fi\fi\fi\fi}%
 \ \FN@\next@}
\def\mdots@{\FN@\mdots@@}
\def\mdots@@{\gdef\thedots@{\dotso@}
 \ifx\next\boldkey\gdef\thedots@\boldkey{\boldkeydots@}\else                
 \ifx\next\boldsymbol\gdef\thedots@\boldsymbol{\boldsymboldots@}\else       
 \ifx,\next\gdef\thedots@{\dotsc}
 \else\ifx\not\next\gdef\thedots@{\dotsb@}
 \else\keybin@
 \ifkeybin@\gdef\thedots@{\dotsb@}
 \else\xdef\meaning@{\meaning\next..........}\xdef\meaning@@{\meaning@}
  \expandafter\math@\meaning@\math@
  \ifmath@
   \expandafter\mathch@\meaning@\mathch@
   \ifmathch@\expandafter\getmathch@\meaning@\getmathch@\fi                 
  \else\expandafter\macro@\meaning@@\macro@                                 
  \ifmacro@                                                                
   \expandafter\not@\meaning@\not@\ifnot@\gdef\thedots@{\dotsb@}
  \else\expandafter\DOTS@\meaning@\DOTS@
  \ifDOTS@
   \ifcase\number\DOTSCASE@\gdef\thedots@{\dotsb@}%
    \or\gdef\thedots@{\dotsi}\else\fi                                      
  \else\expandafter\math@\meaning@\math@                                   
  \ifmath@\expandafter\mathbin@\meaning@\mathbin@
  \ifmathbin@\gdef\thedots@{\dotsb@}
  \else\expandafter\mathrel@\meaning@\mathrel@
  \ifmathrel@\gdef\thedots@{\dotsb@}
  \fi\fi\fi\fi\fi\fi\fi\fi\fi\fi\fi\fi
 \thedots@}
\def\plainldots@{\mathinner{\ldotp\ldotp\ldotp}}
\def\plaincdots@{\mathinner{\cdotp\cdotp\cdotp}}
\def\dotsi{\!\plaincdots@}
\let\dotsb@\plaincdots@
\newif\ifextra@
\newif\ifrightdelim@
\def\rightdelim@{\global\rightdelim@true                                    
 \ifx\next)\else                                                            
 \ifx\next]\else
 \ifx\next\rbrack\else
 \ifx\next\}\else
 \ifx\next\rbrace\else
 \ifx\next\rangle\else
 \ifx\next\rceil\else
 \ifx\next\rfloor\else
 \ifx\next\rgroup\else
 \ifx\next\rmoustache\else
 \ifx\next\right\else
 \ifx\next\bigr\else
 \ifx\next\biggr\else
 \ifx\next\Bigr\else                                                        
 \ifx\next\Biggr\else\global\rightdelim@false
 \fi\fi\fi\fi\fi\fi\fi\fi\fi\fi\fi\fi\fi\fi\fi}
\def\extra@{%
 \global\extra@false\rightdelim@\ifrightdelim@\global\extra@true            
 \else\ifx\next$\global\extra@true                                          
 \else\xdef\meaning@{\meaning\next..........}
 \expandafter\macro@\meaning@\macro@\ifmacro@                               
 \expandafter\DOTS@\meaning@\DOTS@
 \ifDOTS@
 \ifnum\DOTSCASE@=\tw@\global\extra@true                                    
 \fi\fi\fi\fi\fi}
\newif\ifbold@
\def\dotso@{\relaxnext@
 \ifbold@
  \let\next\delayed@
  \DNii@{\extra@\plainldots@\ifextra@\,\fi}%
 \else
  \DNii@{\DN@{\extra@\plainldots@\ifextra@\,\fi}\FN@\next@}%
 \fi
 \nextii@}
\def\extrap@#1{%
 \ifx\next,\DN@{#1\,}\else
 \ifx\next;\DN@{#1\,}\else
 \ifx\next.\DN@{#1\,}\else\extra@
 \ifextra@\DN@{#1\,}\else
 \let\next@#1\fi\fi\fi\fi\next@}
\def\ldots{\DN@{\extrap@\plainldots@}%
 \FN@\next@}
\def\cdots{\DN@{\extrap@\plaincdots@}%
 \FN@\next@}

\def\dotsc{\relaxnext@
 \DN@{\ifx\next;\plainldots@\,\else
  \ifx\next.\plainldots@\,\else\extra@\plainldots@
  \ifextra@\,\fi\fi\fi}%
 \FN@\next@}
\def\cdot{\mathchar"2201 }

\def\mapsto{\DOTSB\mapstochar\rightarrow}

\message{special superscripts,}
\def\dddot#1{{\mathop{#1}\limits^{\vbox to-1.4\ex@{\kern-\tw@\ex@
 \hbox{\rm...}\vss}}}}
\def\ddddot#1{{\mathop{#1}\limits^{\vbox to-1.4\ex@{\kern-\tw@\ex@
 \hbox{\rm....}\vss}}}}
\def\sphat{^{\mathchoice{}{}%
 {\,\,\botsmash{\hbox{\lower4\ex@\hbox{$\m@th\widehat{\null}$}}}}%
 {\,\botsmash{\hbox{\lower3\ex@\hbox{$\m@th\hat{\null}$}}}}}}

\def\spacute{^{\!\botsmash{\hbox{\lower\@ne ex\hbox{\'{}}}}}}
\def\spgrave{^{\mathchoice{}{}{}{\!}%
 \botsmash{\hbox{\lower\@ne ex\hbox{\`{}}}}}}
\def\spdot{^{\hbox{\raise\ex@\hbox{\rm.}}}}
\def\spddot{^{\hbox{\raise\ex@\hbox{\rm..}}}}
\def\spdddot{^{\hbox{\raise\ex@\hbox{\rm...}}}}
\def\spddddot{^{\hbox{\raise\ex@\hbox{\rm....}}}}
\def\spbreve{^{\!\botsmash{\hbox{\lower4\ex@\hbox{\u{}}}}}}

\message{\string\text,}
\def\textonlyfont@#1#2{\def#1{\RIfM@
 \Err@{Use \string#1\space only in text}\else#2\fi}}
\textonlyfont@\rm\tenrm
\textonlyfont@\it\tenit
\textonlyfont@\sl\tensl
\textonlyfont@\bf\tenbf
\def\oldnos#1{\RIfM@{\mathcode`\,="013B \fam\@ne#1}\else
 \leavevmode\hbox{$\m@th\mathcode`\,="013B \fam\@ne#1$}\fi}
\def\text{\RIfM@\expandafter\text@\else\expandafter\text@@\fi}
\def\text@@#1{\leavevmode\hbox{#1}}
\def\mathhexbox@#1#2#3{\text{$\m@th\mathchar"#1#2#3$}}
\def\dag{{\mathhexbox@279}}
\def\ddag{{\mathhexbox@27A}}
\def\S{{\mathhexbox@278}}
\def\P{{\mathhexbox@27B}}
\newif\iffirstchoice@
\firstchoice@true
\def\text@#1{\mathchoice
 {\hbox{\everymath{\displaystyle}\def\textfonti{\the\textfont\@ne}%
  \def\textfontii{\the\textfont\tw@}\textdef@@ T#1}}
 {\hbox{\firstchoice@false
  \everymath{\textstyle}\def\textfonti{\the\textfont\@ne}%
  \def\textfontii{\the\textfont\tw@}\textdef@@ T#1}}
 {\hbox{\firstchoice@false
  \everymath{\scriptstyle}\def\textfonti{\the\scriptfont\@ne}%
  \def\textfontii{\the\scriptfont\tw@}\textdef@@ S\rm#1}}
 {\hbox{\firstchoice@false
  \everymath{\scriptscriptstyle}\def\textfonti
  {\the\scriptscriptfont\@ne}%
  \def\textfontii{\the\scriptscriptfont\tw@}\textdef@@ s\rm#1}}}
\def\textdef@@#1{\textdef@#1\rm\textdef@#1\bf\textdef@#1\sl\textdef@#1\it}
\def\rmfam{0}
\def\textdef@#1#2{%
 \DN@{\csname\expandafter\eat@\string#2fam\endcsname}%
 \if S#1\edef#2{\the\scriptfont\next@\relax}%
 \else\if s#1\edef#2{\the\scriptscriptfont\next@\relax}%
 \else\edef#2{\the\textfont\next@\relax}\fi\fi}
\scriptfont\itfam\tenit \scriptscriptfont\itfam\tenit
\scriptfont\slfam\tensl \scriptscriptfont\slfam\tensl
\newif\iftopfolded@
\newif\ifbotfolded@
\def\topfoldedtext{\topfolded@true\botfolded@false\foldedtext@}
\def\botfoldedtext{\botfolded@true\topfolded@false\foldedtext@}
\def\foldedtext{\topfolded@false\botfolded@false\foldedtext@}
\Invalid@\foldedwidth
\def\foldedtext@{\relaxnext@
 \DN@{\ifx\next\foldedwidth\let\next@\nextii@\else
  \DN@{\nextii@\foldedwidth{.3\hsize}}\fi\next@}%
 \DNii@\foldedwidth##1##2{\setbox\z@\vbox
  {\normalbaselines\hsize##1\relax
  \tolerance1600 \noindent\ignorespaces##2}\ifbotfolded@\boxz@\else
  \iftopfolded@\vtop{\unvbox\z@}\else\vcenter{\boxz@}\fi\fi}%
 \FN@\next@}
\message{math font commands,}
\def\bold{\RIfM@\expandafter\bold@\else
 \expandafter\nonmatherr@\expandafter\bold\fi}
\def\bold@#1{{\bold@@{#1}}}
\def\bold@@#1{\fam\bffam\relax#1}
\def\slanted{\RIfM@\expandafter\slanted@\else
 \expandafter\nonmatherr@\expandafter\slanted\fi}
\def\slanted@#1{{\slanted@@{#1}}}
\def\slanted@@#1{\fam\slfam\relax#1}
\def\roman{\RIfM@\expandafter\roman@\else
 \expandafter\nonmatherr@\expandafter\roman\fi}
\def\roman@#1{{\roman@@{#1}}}
\def\roman@@#1{\fam\rmfam\relax#1}
\def\italic{\RIfM@\expandafter\italic@\else
 \expandafter\nonmatherr@\expandafter\italic\fi}
\def\italic@#1{{\italic@@{#1}}}
\def\italic@@#1{\fam\itfam\relax#1}
\def\Cal{\RIfM@\expandafter\Cal@\else
 \expandafter\nonmatherr@\expandafter\Cal\fi}
\def\Cal@#1{{\Cal@@{#1}}}
\def\Cal@@#1{\noaccents@\fam\tw@#1}
\mathchardef\Gamma="0000
\mathchardef\Delta="0001
\mathchardef\Theta="0002
\mathchardef\Lambda="0003
\mathchardef\Xi="0004
\mathchardef\Pi="0005
\mathchardef\Sigma="0006
\mathchardef\Upsilon="0007
\mathchardef\Phi="0008
\mathchardef\Psi="0009
\mathchardef\Omega="000A
\mathchardef\varGamma="0100
\mathchardef\varDelta="0101
\mathchardef\varTheta="0102
\mathchardef\varLambda="0103
\mathchardef\varXi="0104
\mathchardef\varPi="0105
\mathchardef\varSigma="0106
\mathchardef\varUpsilon="0107
\mathchardef\varPhi="0108
\mathchardef\varPsi="0109
\mathchardef\varOmega="010A
\let\alloc@@\alloc@
\def\hexnumber@#1{\ifcase#1 0\or 1\or 2\or 3\or 4\or 5\or 6\or 7\or 8\or
 9\or A\or B\or C\or D\or E\or F\fi}
\def\loadmsam{%
 \font@\tenmsa=msam10
 \font@\sevenmsa=msam7
 \font@\fivemsa=msam5
 \alloc@@8\fam\chardef\sixt@@n\msafam
 \textfont\msafam=\tenmsa
 \scriptfont\msafam=\sevenmsa
 \scriptscriptfont\msafam=\fivemsa
 \edef\next{\hexnumber@\msafam}%
 \mathchardef\dabar@"0\next39
 \edef\dashrightarrow{\mathrel{\dabar@\dabar@\mathchar"0\next4B}}%
 \edef\dashleftarrow{\mathrel{\mathchar"0\next4C\dabar@\dabar@}}%
 \let\dasharrow\dashrightarrow
 \edef\ulcorner{\delimiter"4\next70\next70 }%
 \edef\urcorner{\delimiter"5\next71\next71 }%
 \edef\llcorner{\delimiter"4\next78\next78 }%
 \edef\lrcorner{\delimiter"5\next79\next79 }%
 \edef\yen{{\noexpand\mathhexbox@\next55}}%
 \edef\checkmark{{\noexpand\mathhexbox@\next58}}%
 \edef\circledR{{\noexpand\mathhexbox@\next72}}%
 \edef\maltese{{\noexpand\mathhexbox@\next7A}}%
 \global\let\loadmsam\empty}%
\def\loadmsbm{%
 \font@\tenmsb=msbm10 \font@\sevenmsb=msbm7 \font@\fivemsb=msbm5
 \alloc@@8\fam\chardef\sixt@@n\msbfam
 \textfont\msbfam=\tenmsb
 \scriptfont\msbfam=\sevenmsb \scriptscriptfont\msbfam=\fivemsb
 \global\let\loadmsbm\empty
 }
\def\widehat#1{\ifx\undefined\msbfam \DN@{362}%
  \else \setboxz@h{$\m@th#1$}%
    \edef\next@{\ifdim\wdz@>\tw@ em%
        \hexnumber@\msbfam 5B%
      \else 362\fi}\fi
  \mathaccent"0\next@{#1}}
\def\widetilde#1{\ifx\undefined\msbfam \DN@{365}%
  \else \setboxz@h{$\m@th#1$}%
    \edef\next@{\ifdim\wdz@>\tw@ em%
        \hexnumber@\msbfam 5D%
      \else 365\fi}\fi
  \mathaccent"0\next@{#1}}
\message{\string\newsymbol,}
\def\newsymbol#1#2#3#4#5{\define#1{}%
  \count@#2\relax \advance\count@\m@ne 
 \ifcase\count@
   \ifx\undefined\msafam\loadmsam\fi \let\next@\msafam
 \or \ifx\undefined\msbfam\loadmsbm\fi \let\next@\msbfam
 \else  \Err@{\Invalid@@\string\newsymbol}\let\next@\tw@\fi
 \mathchardef#1="#3\hexnumber@\next@#4#5\space}

\def\Bbb{\RIfM@\expandafter\Bbb@\else
 \expandafter\nonmatherr@\expandafter\Bbb\fi}
\def\Bbb@#1{{\Bbb@@{#1}}}
\def\Bbb@@#1{\noaccents@\fam\msbfam\relax#1}
\message{bold Greek and bold symbols,}
\def\loadbold{%
 \font@\tencmmib=cmmib10 \font@\sevencmmib=cmmib7 \font@\fivecmmib=cmmib5
 \skewchar\tencmmib'177 \skewchar\sevencmmib'177 \skewchar\fivecmmib'177
 \alloc@@8\fam\chardef\sixt@@n\cmmibfam
 \textfont\cmmibfam\tencmmib
 \scriptfont\cmmibfam\sevencmmib \scriptscriptfont\cmmibfam\fivecmmib
 \font@\tencmbsy=cmbsy10 \font@\sevencmbsy=cmbsy7 \font@\fivecmbsy=cmbsy5
 \skewchar\tencmbsy'60 \skewchar\sevencmbsy'60 \skewchar\fivecmbsy'60
 \alloc@@8\fam\chardef\sixt@@n\cmbsyfam
 \textfont\cmbsyfam\tencmbsy
 \scriptfont\cmbsyfam\sevencmbsy \scriptscriptfont\cmbsyfam\fivecmbsy
 \let\loadbold\empty
}
\def\boldnotloaded#1{\Err@{\ifcase#1\or First\else Second\fi
       bold symbol font not loaded}}
\def\mathchari@#1#2#3{\ifx\undefined\cmmibfam
    \boldnotloaded@\@ne
  \else\mathchar"#1\hexnumber@\cmmibfam#2#3\space \fi}
\def\mathcharii@#1#2#3{\ifx\undefined\cmbsyfam
    \boldnotloaded\tw@
  \else \mathchar"#1\hexnumber@\cmbsyfam#2#3\space\fi}
\edef\bffam@{\hexnumber@\bffam}
\def\boldkey#1{\ifcat\noexpand#1A%
  \ifx\undefined\cmmibfam \boldnotloaded\@ne
  \else {\fam\cmmibfam#1}\fi
 \else
 \ifx#1!\mathchar"5\bffam@21 \else
 \ifx#1(\mathchar"4\bffam@28 \else\ifx#1)\mathchar"5\bffam@29 \else
 \ifx#1+\mathchar"2\bffam@2B \else\ifx#1:\mathchar"3\bffam@3A \else
 \ifx#1;\mathchar"6\bffam@3B \else\ifx#1=\mathchar"3\bffam@3D \else
 \ifx#1?\mathchar"5\bffam@3F \else\ifx#1[\mathchar"4\bffam@5B \else
 \ifx#1]\mathchar"5\bffam@5D \else
 \ifx#1,\mathchari@63B \else
 \ifx#1-\mathcharii@200 \else
 \ifx#1.\mathchari@03A \else
 \ifx#1/\mathchari@03D \else
 \ifx#1<\mathchari@33C \else
 \ifx#1>\mathchari@33E \else
 \ifx#1*\mathcharii@203 \else
 \ifx#1|\mathcharii@06A \else
 \ifx#10\bold0\else\ifx#11\bold1\else\ifx#12\bold2\else\ifx#13\bold3\else
 \ifx#14\bold4\else\ifx#15\bold5\else\ifx#16\bold6\else\ifx#17\bold7\else
 \ifx#18\bold8\else\ifx#19\bold9\else
  \Err@{\string\boldkey\space can't be used with #1}%
 \fi\fi\fi\fi\fi\fi\fi\fi\fi\fi\fi\fi\fi\fi\fi
 \fi\fi\fi\fi\fi\fi\fi\fi\fi\fi\fi\fi\fi\fi}
\def\boldsymbol#1{%
 \DN@{\Err@{You can't use \string\boldsymbol\space with \string#1}#1}%
 \ifcat\noexpand#1A%
   \let\next@\relax
   \ifx\undefined\cmmibfam \boldnotloaded\@ne
   \else {\fam\cmmibfam#1}\fi
 \else
  \xdef\meaning@{\meaning#1.........}%
  \expandafter\math@\meaning@\math@
  \ifmath@
   \expandafter\mathch@\meaning@\mathch@
   \ifmathch@
    \expandafter\boldsymbol@@\meaning@\boldsymbol@@
   \fi
  \else
   \expandafter\macro@\meaning@\macro@
   \expandafter\delim@\meaning@\delim@
   \ifdelim@
    \expandafter\delim@@\meaning@\delim@@
   \else
    \boldsymbol@{#1}%
   \fi
  \fi
 \fi
 \next@}
\def\mathhexboxii@#1#2{\ifx\undefined\cmbsyfam
    \boldnotloaded\tw@
  \else \mathhexbox@{\hexnumber@\cmbsyfam}{#1}{#2}\fi}
\def\boldsymbol@#1{\let\next@\relax\let\next#1%
 \ifx\next\cdot\mathcharii@201 \else
 \ifx\next\prime{{\null\mathcharii@030 \null}}\else
 \ifx\next\lbrack\mathchar"4\bffam@5B \else
 \ifx\next\rbrack\mathchar"5\bffam@5D \else
 \ifx\next\{\mathcharii@466 \else
 \ifx\next\lbrace\mathcharii@466 \else
 \ifx\next\}\mathcharii@567 \else
 \ifx\next\rbrace\mathcharii@567 \else
 \ifx\next\surd{{\mathcharii@170}}\else
 \ifx\next\S{{\mathhexboxii@78}}\else
 \ifx\next\P{{\mathhexboxii@7B}}\else
 \ifx\next\dag{{\mathhexboxii@79}}\else
 \ifx\next\ddag{{\mathhexboxii@7A}}\else
 \DN@{\Err@{You can't use \string\boldsymbol\space with \string#1}#1}%
 \fi\fi\fi\fi\fi\fi\fi\fi\fi\fi\fi\fi\fi}
\def\boldsymbol@@#1.#2\boldsymbol@@{\classnum@#1 \count@@@\classnum@        
 \divide\classnum@4096 \count@\classnum@                                    
 \multiply\count@4096 \advance\count@@@-\count@ \count@@\count@@@           
 \divide\count@@@\@cclvi \count@\count@@                                    
 \multiply\count@@@\@cclvi \advance\count@@-\count@@@                       
 \divide\count@@@\@cclvi                                                    
 \multiply\classnum@4096 \advance\classnum@\count@@                         
 \ifnum\count@@@=\z@                                                        
  \count@"\bffam@ \multiply\count@\@cclvi
  \advance\classnum@\count@
  \DN@{\mathchar\number\classnum@}%
 \else
  \ifnum\count@@@=\@ne                                                      
   \ifx\undefined\cmmibfam \DN@{\boldnotloaded\@ne}%
   \else \count@\cmmibfam \multiply\count@\@cclvi
     \advance\classnum@\count@
     \DN@{\mathchar\number\classnum@}\fi
  \else
   \ifnum\count@@@=\tw@                                                    
     \ifx\undefined\cmbsyfam
       \DN@{\boldnotloaded\tw@}%
     \else
       \count@\cmbsyfam \multiply\count@\@cclvi
       \advance\classnum@\count@
       \DN@{\mathchar\number\classnum@}%
     \fi
  \fi
 \fi
\fi}
\newif\ifdelim@
\newcount\delimcount@
{\uccode`6=`\\ \uccode`7=`d \uccode`8=`e \uccode`9=`l
 \uppercase{\gdef\delim@#1#2#3#4#5\delim@
  {\delim@false\ifx 6#1\ifx 7#2\ifx 8#3\ifx 9#4\delim@true
   \xdef\meaning@{#5}\fi\fi\fi\fi}}}
\def\delim@@#1"#2#3#4#5#6\delim@@{\if#32%
\let\next@\relax
 \ifx\undefined\cmbsyfam \boldnotloaded\@ne
 \else \mathcharii@#2#4#5\space \fi\fi}
\def\vert{\delimiter"026A30C }
\def\Vert{\delimiter"026B30D }
\let\|\Vert
\def\backslash{\delimiter"026E30F }
\def\boldkeydots@#1{\bold@true\let\next=#1\let\delayed@=#1\mdots@@
 \boldkey#1\bold@false}  
\def\boldsymboldots@#1{\bold@true\let\next#1\let\delayed@#1\mdots@@
 \boldsymbol#1\bold@false}
\message{Euler fonts,}

\def\frak{\mathfont@\frak}

\def\loadmathfont#1{%
   \expandafter\font@\csname ten#1\endcsname=#110
   \expandafter\font@\csname seven#1\endcsname=#17
   \expandafter\font@\csname five#1\endcsname=#15
   \edef\next{\noexpand\alloc@@8\fam\chardef\sixt@@n
     \expandafter\noexpand\csname#1fam\endcsname}%
   \next
   \textfont\csname#1fam\endcsname \csname ten#1\endcsname
   \scriptfont\csname#1fam\endcsname \csname seven#1\endcsname
   \scriptscriptfont\csname#1fam\endcsname \csname five#1\endcsname
   \expandafter\def\csname #1\expandafter\endcsname\expandafter{%
      \expandafter\mathfont@\csname#1\endcsname}%
 \expandafter\gdef\csname load#1\endcsname{}%
}
\def\mathfont@#1{\RIfM@\expandafter\mathfont@@\expandafter#1\else
  \expandafter\nonmatherr@\expandafter#1\fi}
\def\mathfont@@#1#2{{\mathfont@@@#1{#2}}}
\def\mathfont@@@#1#2{\noaccents@
   \fam\csname\expandafter\eat@\string#1fam\endcsname
   \relax#2}
\message{math accents,}
\def\accentclass@{7}
\def\noaccents@{\def\accentclass@{0}}
\def\makeacc@#1#2{\def#1{\mathaccent"\accentclass@#2 }}
\makeacc@\hat{05E}
\makeacc@\check{014}
\makeacc@\tilde{07E}
\makeacc@\acute{013}
\makeacc@\grave{012}
\makeacc@\dot{05F}
\makeacc@\ddot{07F}
\makeacc@\breve{015}
\makeacc@\bar{016}

\newcount\skewcharcount@
\newcount\familycount@
\def\theskewchar@{\familycount@\@ne
 \global\skewcharcount@\the\skewchar\textfont\@ne                           
 \ifnum\fam>\m@ne\ifnum\fam<16
  \global\familycount@\the\fam\relax
  \global\skewcharcount@\the\skewchar\textfont\the\fam\relax\fi\fi          
 \ifnum\skewcharcount@>\m@ne
  \ifnum\skewcharcount@<128
  \multiply\familycount@256
  \global\advance\skewcharcount@\familycount@
  \global\advance\skewcharcount@28672
  \mathchar\skewcharcount@\else
  \global\skewcharcount@\m@ne\fi\else
 \global\skewcharcount@\m@ne\fi}                                            
\newcount\pointcount@
\def\getpoints@#1.#2\getpoints@{\pointcount@#1 }
\newdimen\accentdimen@
\newcount\accentmu@
\def\dimentomu@{\multiply\accentdimen@ 100
 \expandafter\getpoints@\the\accentdimen@\getpoints@
 \multiply\pointcount@18
 \divide\pointcount@\@m
 \global\accentmu@\pointcount@}
\def\Makeacc@#1#2{\def#1{\RIfM@\DN@{\mathaccent@
 {"\accentclass@#2 }}\else\DN@{\nonmatherr@{#1}}\fi\next@}}
\def\unbracefonts@{\let\Cal@\Cal@@\let\roman@\roman@@\let\bold@\bold@@
 \let\slanted@\slanted@@}
\def\mathaccent@#1#2{\ifnum\fam=\m@ne\xdef\thefam@{1}\else
 \xdef\thefam@{\the\fam}\fi                                                 
 \accentdimen@\z@                                                           
 \setboxz@h{\unbracefonts@$\m@th\fam\thefam@\relax#2$}
 \ifdim\accentdimen@=\z@\DN@{\mathaccent#1{#2}}
  \setbox@ne\hbox{\unbracefonts@$\m@th\fam\thefam@\relax#2\theskewchar@$}
  \setbox\tw@\hbox{$\m@th\ifnum\skewcharcount@=\m@ne\else
   \mathchar\skewcharcount@\fi$}
  \global\accentdimen@\wd@ne\global\advance\accentdimen@-\wdz@
  \global\advance\accentdimen@-\wd\tw@                                     
  \global\multiply\accentdimen@\tw@
  \dimentomu@\global\advance\accentmu@\@ne                                 
 \else\DN@{{\mathaccent#1{#2\mkern\accentmu@ mu}%
    \mkern-\accentmu@ mu}{}}\fi                                             
 \next@}\Makeacc@\Hat{05E}
\Makeacc@\Check{014}
\Makeacc@\Tilde{07E}
\Makeacc@\Acute{013}
\Makeacc@\Grave{012}
\Makeacc@\Dot{05F}
\Makeacc@\Ddot{07F}
\Makeacc@\Breve{015}
\Makeacc@\Bar{016}
\def\Vec{\RIfM@\DN@{\mathaccent@{"017E }}\else
 \DN@{\nonmatherr@\Vec}\fi\next@}
\def\accentedsymbol#1#2{\csname newbox\expandafter\endcsname
  \csname\expandafter\eat@\string#1@box\endcsname
 \expandafter\setbox\csname\expandafter\eat@
  \string#1@box\endcsname\hbox{$\m@th#2$}\define
  #1{\copy\csname\expandafter\eat@\string#1@box\endcsname{}}}
\message{roots,}
\def\sqrt#1{\radical"270370 {#1}}
\let\underline@\underline
\let\overline@\overline
\def\underline#1{\underline@{#1}}
\def\overline#1{\overline@{#1}}
\Invalid@\leftroot
\Invalid@\uproot
\newcount\uproot@
\newcount\leftroot@
\def\root{\relaxnext@
  \DN@{\ifx\next\uproot\let\next@\nextii@\else
   \ifx\next\leftroot\let\next@\nextiii@\else
   \let\next@\plainroot@\fi\fi\next@}%
  \DNii@\uproot##1{\uproot@##1\relax\FN@\nextiv@}%
  \def\nextiv@{\ifx\next\space@\DN@. {\FN@\nextv@}\else
   \DN@.{\FN@\nextv@}\fi\next@.}%
  \def\nextv@{\ifx\next\leftroot\let\next@\nextvi@\else
   \let\next@\plainroot@\fi\next@}%
  \def\nextvi@\leftroot##1{\leftroot@##1\relax\plainroot@}%
   \def\nextiii@\leftroot##1{\leftroot@##1\relax\FN@\nextvii@}%
  \def\nextvii@{\ifx\next\space@
   \DN@. {\FN@\nextviii@}\else
   \DN@.{\FN@\nextviii@}\fi\next@.}%
  \def\nextviii@{\ifx\next\uproot\let\next@\nextix@\else
   \let\next@\plainroot@\fi\next@}%
  \def\nextix@\uproot##1{\uproot@##1\relax\plainroot@}%
  \bgroup\uproot@\z@\leftroot@\z@\FN@\next@}
\def\plainroot@#1\of#2{\setbox\rootbox\hbox{$\m@th\scriptscriptstyle{#1}$}%
 \mathchoice{\r@@t\displaystyle{#2}}{\r@@t\textstyle{#2}}
 {\r@@t\scriptstyle{#2}}{\r@@t\scriptscriptstyle{#2}}\egroup}
\def\r@@t#1#2{\setboxz@h{$\m@th#1\sqrt{#2}$}%
 \dimen@\ht\z@\advance\dimen@-\dp\z@
 \setbox@ne\hbox{$\m@th#1\mskip\uproot@ mu$}\advance\dimen@ 1.667\wd@ne
 \mkern-\leftroot@ mu\mkern5mu\raise.6\dimen@\copy\rootbox
 \mkern-10mu\mkern\leftroot@ mu\boxz@}
\def\boxed#1{\setboxz@h{$\m@th\displaystyle{#1}$}\dimen@.4\ex@
 \advance\dimen@3\ex@\advance\dimen@\dp\z@
 \hbox{\lower\dimen@\hbox{%
 \vbox{\hrule height.4\ex@
 \hbox{\vrule width.4\ex@\hskip3\ex@\vbox{\vskip3\ex@\boxz@\vskip3\ex@}%
 \hskip3\ex@\vrule width.4\ex@}\hrule height.4\ex@}%
 }}}
\message{commutative diagrams,}
\let\ampersand@\relax
\newdimen\minaw@
\minaw@11.11128\ex@
\newdimen\minCDaw@
\minCDaw@2.5pc
\def\minCDarrowwidth#1{\RIfMIfI@\onlydmatherr@\minCDarrowwidth
 \else\minCDaw@#1\relax\fi\else\onlydmatherr@\minCDarrowwidth\fi}
\newif\ifCD@
\def\CD{\bgroup\vspace@\relax\let\ampersand@&\iffalse}\fi
 \CD@true\vcenter\bgroup\Let@\tabskip\z@skip\baselineskip20\ex@
 \lineskip3\ex@\lineskiplimit3\ex@\halign\bgroup
 &\hfill$\m@th##$\hfill\crcr}
\def\endCD{\crcr\egroup\egroup\egroup}
\newdimen\bigaw@
\atdef@>#1>#2>{\ampersand@                                                  
 \setboxz@h{$\m@th\ssize\;{#1}\;\;$}
 \setbox@ne\hbox{$\m@th\ssize\;{#2}\;\;$}
 \setbox\tw@\hbox{$\m@th#2$}
 \ifCD@\global\bigaw@\minCDaw@\else\global\bigaw@\minaw@\fi                 
 \ifdim\wdz@>\bigaw@\global\bigaw@\wdz@\fi
 \ifdim\wd@ne>\bigaw@\global\bigaw@\wd@ne\fi                                
 \ifCD@\enskip\fi                                                           
 \ifdim\wd\tw@>\z@
  \mathrel{\mathop{\hbox to\bigaw@{\rightarrowfill@\displaystyle}}%
    \limits^{#1}_{#2}}
 \else\mathrel{\mathop{\hbox to\bigaw@{\rightarrowfill@\displaystyle}}%
    \limits^{#1}}\fi                                                        
 \ifCD@\enskip\fi                                                          
 \ampersand@}                                                              
\atdef@<#1<#2<{\ampersand@\setboxz@h{$\m@th\ssize\;\;{#1}\;$}%
 \setbox@ne\hbox{$\m@th\ssize\;\;{#2}\;$}\setbox\tw@\hbox{$\m@th#2$}%
 \ifCD@\global\bigaw@\minCDaw@\else\global\bigaw@\minaw@\fi
 \ifdim\wdz@>\bigaw@\global\bigaw@\wdz@\fi
 \ifdim\wd@ne>\bigaw@\global\bigaw@\wd@ne\fi
 \ifCD@\enskip\fi
 \ifdim\wd\tw@>\z@
  \mathrel{\mathop{\hbox to\bigaw@{\leftarrowfill@\displaystyle}}%
       \limits^{#1}_{#2}}\else
  \mathrel{\mathop{\hbox to\bigaw@{\leftarrowfill@\displaystyle}}%
       \limits^{#1}}\fi
 \ifCD@\enskip\fi\ampersand@}
\begingroup
 \catcode`\~=\active \lccode`\~=`\@
 \lowercase{%
  \global\atdef@)#1)#2){~>#1>#2>}
  \global\atdef@(#1(#2({~<#1<#2<}}
\endgroup
\atdef@ A#1A#2A{\llap{$\m@th\vcenter{\hbox
 {$\ssize#1$}}$}\Big\uparrow\rlap{$\m@th\vcenter{\hbox{$\ssize#2$}}$}&&}
\atdef@ V#1V#2V{\llap{$\m@th\vcenter{\hbox
 {$\ssize#1$}}$}\Big\downarrow\rlap{$\m@th\vcenter{\hbox{$\ssize#2$}}$}&&}
\atdef@={&\enskip\mathrel
 {\vbox{\hrule width\minCDaw@\vskip3\ex@\hrule width
 \minCDaw@}}\enskip&}
\atdef@|{\Big\Vert&&}
\atdef@\vert{\Big\Vert&&}
\def\pretend#1\haswidth#2{\setboxz@h{$\m@th\scriptstyle{#2}$}\hbox
 to\wdz@{\hfill$\m@th\scriptstyle{#1}$\hfill}}
\message{poor man's bold,}
\def\pmb{\RIfM@\expandafter\mathpalette\expandafter\pmb@\else
 \expandafter\pmb@@\fi}
\def\pmb@@#1{\leavevmode\setboxz@h{#1}%
   \dimen@-\wdz@
   \kern-.5\ex@\copy\z@
   \kern\dimen@\kern.25\ex@\raise.4\ex@\copy\z@
   \kern\dimen@\kern.25\ex@\box\z@
}
\def\binrel@@#1{\ifdim\wd2<\z@\mathbin{#1}\else\ifdim\wd\tw@>\z@
 \mathrel{#1}\else{#1}\fi\fi}
\newdimen\pmbraise@
\def\pmb@#1#2{\setbox\thr@@\hbox{$\m@th#1{#2}$}%
 \setbox4\hbox{$\m@th#1\mkern.5mu$}\pmbraise@\wd4\relax
 \binrel@{#2}%
 \dimen@-\wd\thr@@
   \binrel@@{%
   \mkern-.8mu\copy\thr@@
   \kern\dimen@\mkern.4mu\raise\pmbraise@\copy\thr@@
   \kern\dimen@\mkern.4mu\box\thr@@
}}
\def\documentstyle#1{\W@{}\input #1.sty\relax}
\message{syntax check,}
\font\dummyft@=dummy
\fontdimen1 \dummyft@=\z@
\fontdimen2 \dummyft@=\z@
\fontdimen3 \dummyft@=\z@
\fontdimen4 \dummyft@=\z@
\fontdimen5 \dummyft@=\z@
\fontdimen6 \dummyft@=\z@
\fontdimen7 \dummyft@=\z@
\fontdimen8 \dummyft@=\z@
\fontdimen9 \dummyft@=\z@
\fontdimen10 \dummyft@=\z@
\fontdimen11 \dummyft@=\z@
\fontdimen12 \dummyft@=\z@
\fontdimen13 \dummyft@=\z@
\fontdimen14 \dummyft@=\z@
\fontdimen15 \dummyft@=\z@
\fontdimen16 \dummyft@=\z@
\fontdimen17 \dummyft@=\z@
\fontdimen18 \dummyft@=\z@
\fontdimen19 \dummyft@=\z@
\fontdimen20 \dummyft@=\z@
\fontdimen21 \dummyft@=\z@
\fontdimen22 \dummyft@=\z@
\def\fontlist@{\\{\tenrm}\\{\sevenrm}\\{\fiverm}\\{\teni}\\{\seveni}%
 \\{\fivei}\\{\tensy}\\{\sevensy}\\{\fivesy}\\{\tenex}\\{\tenbf}\\{\sevenbf}%
 \\{\fivebf}\\{\tensl}\\{\tenit}}
\def\font@#1=#2 {\rightappend@#1\to\fontlist@\font#1=#2 }
\def\dodummy@{{\def\\##1{\global\let##1\dummyft@}\fontlist@}}
\def\nopages@{\output{\setbox\z@\box\@cclv \deadcycles\z@}%
 \alloc@5\toks\toksdef\@cclvi\output}
\let\galleys\nopages@
\newif\ifsyntax@
\newcount\countxviii@
\def\syntax{\syntax@true\dodummy@\countxviii@\count18
 \loop\ifnum\countxviii@>\m@ne\textfont\countxviii@=\dummyft@
 \scriptfont\countxviii@=\dummyft@\scriptscriptfont\countxviii@=\dummyft@
 \advance\countxviii@\m@ne\repeat                                           
 \dummyft@\tracinglostchars\z@\nopages@\frenchspacing\hbadness\@M}
\def\first@#1#2\end{#1}
\def\printoptions{\W@{Do you want S(yntax check),
  G(alleys) or P(ages)?}%
 \message{Type S, G or P, followed by <return>: }%
 \begingroup 
 \endlinechar\m@ne 
 \read\m@ne to\ans@
 \edef\ans@{\uppercase{\def\noexpand\ans@{%
   \expandafter\first@\ans@ P\end}}}%
 \expandafter\endgroup\ans@
 \if\ans@ P
 \else \if\ans@ S\syntax
 \else \if\ans@ G\galleys
 \else\message{? Unknown option: \ans@; using the `pages' option.}%
 \fi\fi\fi}
\def\alloc@#1#2#3#4#5{\global\advance\count1#1by\@ne
 \ch@ck#1#4#2\allocationnumber=\count1#1
 \global#3#5=\allocationnumber
 \ifalloc@\wlog{\string#5=\string#2\the\allocationnumber}\fi}
\def\document{\def\alloclist@{}\def\fontlist@{}}
\let\enddocument\bye

\let\proclaim\undefined
\let\footnote\undefined
\let\=\undefined
\let\>\undefined

\catcode`\@=\active
\message{... finished}

\expandafter\ifx\csname mathdefs.tex\endcsname\relax
  \expandafter\gdef\csname mathdefs.tex\endcsname{}
\else \message{Hey!  Apparently you were trying to
  \string\input{mathdefs.tex} twice.   This does not make sense.} 
\errmessage{Please edit your file (probably \jobname.tex) and remove
any duplicate ``\string\input'' lines}\endinput\fi




\catcode`\X=12\catcode`\@=11

\def\n@wcount{\alloc@0\count\countdef\insc@unt}
\def\n@wwrite{\alloc@7\write\chardef\sixt@@n}
\def\n@wread{\alloc@6\read\chardef\sixt@@n}
\def\r@s@t{\relax}\def\v@idline{\par}\def\@mputate#1/{#1}
\def\l@c@l#1X{\firstpart.#1}\def\gl@b@l#1X{#1}\def\t@d@l#1X{{}}

\def\crossrefs#1{\ifx\all#1\let\tr@ce=\all\else\def\tr@ce{#1,}\fi
   \n@wwrite\cit@tionsout\openout\cit@tionsout=\jobname.cit 
   \write\cit@tionsout{\tr@ce}\expandafter\setfl@gs\tr@ce,}
\def\setfl@gs#1,{\def\@{#1}\ifx\@\empty\let\next=\relax
   \else\let\next=\setfl@gs\expandafter\xdef
   \csname#1tr@cetrue\endcsname{}\fi\next}
\def\m@ketag#1#2{\expandafter\n@wcount\csname#2tagno\endcsname
     \csname#2tagno\endcsname=0\let\tail=\all\xdef\all{\tail#2,}
   \ifx#1\l@c@l\let\tail=\r@s@t\xdef\r@s@t{\csname#2tagno\endcsname=0\tail}\fi
   \expandafter\gdef\csname#2cite\endcsname##1{\expandafter
     \ifx\csname#2tag##1\endcsname\relax?\else\csname#2tag##1\endcsname\fi
     \expandafter\ifx\csname#2tr@cetrue\endcsname\relax\else
     \write\cit@tionsout{#2tag ##1 cited on page \folio.}\fi}
   \expandafter\gdef\csname#2page\endcsname##1{\expandafter
     \ifx\csname#2page##1\endcsname\relax?\else\csname#2page##1\endcsname\fi
     \expandafter\ifx\csname#2tr@cetrue\endcsname\relax\else
     \write\cit@tionsout{#2tag ##1 cited on page \folio.}\fi}
   \expandafter\gdef\csname#2tag\endcsname##1{\expandafter
      \ifx\csname#2check##1\endcsname\relax
      \expandafter\xdef\csname#2check##1\endcsname{}%
      \else\immediate\write16{Warning: #2tag ##1 used more than once.}\fi
      \multit@g{#1}{#2}##1/X%
      \write\t@gsout{#2tag ##1 assigned number \csname#2tag##1\endcsname\space
      on page \number\count0.}%
   \csname#2tag##1\endcsname}}

\def\multit@g#1#2#3/#4X{\def\t@mp{#4}\ifx\t@mp\empty%
      \global\advance\csname#2tagno\endcsname by 1 
      \expandafter\xdef\csname#2tag#3\endcsname
      {#1\number\csname#2tagno\endcsnameX}%
   \else\expandafter\ifx\csname#2last#3\endcsname\relax
      \expandafter\n@wcount\csname#2last#3\endcsname
      \global\advance\csname#2tagno\endcsname by 1 
      \expandafter\xdef\csname#2tag#3\endcsname
      {#1\number\csname#2tagno\endcsnameX}
      \write\t@gsout{#2tag #3 assigned number \csname#2tag#3\endcsname\space
      on page \number\count0.}\fi
   \global\advance\csname#2last#3\endcsname by 1
   \def\t@mp{\expandafter\xdef\csname#2tag#3/}%
   \expandafter\t@mp\@mputate#4\endcsname
   {\csname#2tag#3\endcsname\lastpart{\csname#2last#3\endcsname}}\fi}
\def\t@gs#1{\def\all{}\m@ketag#1e\m@ketag#1s\m@ketag\t@d@l p
\let\realscite\scite
\let\realstag\stag
   \m@ketag\gl@b@l r \n@wread\t@gsin
   \openin\t@gsin=\jobname.tgs \re@der \closein\t@gsin
   \n@wwrite\t@gsout\openout\t@gsout=\jobname.tgs }
\outer\def\localtags{\t@gs\l@c@l}
\outer\def\globaltags{\t@gs\gl@b@l}
\outer\def\newlocaltag#1{\m@ketag\l@c@l{#1}}
\outer\def\newglobaltag#1{\m@ketag\gl@b@l{#1}}

\newif\ifpr@ 
\def\m@kecs #1tag #2 assigned number #3 on page #4.%
   {\expandafter\gdef\csname#1tag#2\endcsname{#3}
   \expandafter\gdef\csname#1page#2\endcsname{#4}
   \ifpr@\expandafter\xdef\csname#1check#2\endcsname{}\fi}
\def\re@der{\ifeof\t@gsin\let\next=\relax\else
   \read\t@gsin to\t@gline\ifx\t@gline\v@idline\else
   \expandafter\m@kecs \t@gline\fi\let \next=\re@der\fi\next}
\def\pretags#1{\pr@true\pret@gs#1,,}
\def\pret@gs#1,{\def\@{#1}\ifx\@\empty\let\n@xtfile=\relax
   \else\let\n@xtfile=\pret@gs \openin\t@gsin=#1.tgs \message{#1} \re@der 
   \closein\t@gsin\fi \n@xtfile}

\newcount\sectno\sectno=0\newcount\subsectno\subsectno=0
\newif\ifultr@local \def\ultralocal{\ultr@localtrue}
\def\firstpart{\number\sectno}
\def\lastpart#1{\ifcase#1 \or a\or b\or c\or d\or e\or f\or g\or h\or 
   i\or k\or l\or m\or n\or o\or p\or q\or r\or s\or t\or u\or v\or w\or 
   x\or y\or z \fi}

\def\resetall{\global\advance\sectno by 1\subsectno=0
   \gdef\firstpart{\number\sectno}\r@s@t}
\def\resetsub{\global\advance\subsectno by 1
   \gdef\firstpart{\number\sectno.\number\subsectno}\r@s@t}
\def\newsection#1\par{\resetall\vskip0pt plus.3\vsize\penalty-250
   \vskip0pt plus-.3\vsize\bigskip\bigskip
   \message{#1}\leftline{\bf#1}\nobreak\bigskip}
\def\subsection#1\par{\ifultr@local\resetsub\fi
   \vskip0pt plus.2\vsize\penalty-250\vskip0pt plus-.2\vsize
   \bigskip\smallskip\message{#1}\leftline{\bf#1}\nobreak\medskip}


\newdimen\marginshift

\newdimen\margindelta
\newdimen\marginmax
\newdimen\marginmin

\def\margininit{       
\marginmax=3 true cm                  
				      
\margindelta=0.1 true cm              
\marginmin=0.1true cm                 
\marginshift=\marginmin
}    

\def\t@gsjj#1,{\def\@{#1}\ifx\@\empty\let\next=\relax\else\let\next=\t@gsjj
   \def\@@{p}\ifx\@\@@\else
   \expandafter\gdef\csname#1cite\endcsname##1{\citejj{##1}}
   \expandafter\gdef\csname#1page\endcsname##1{?}
   \expandafter\gdef\csname#1tag\endcsname##1{\tagjj{##1}}\fi\fi\next}
\newif\ifshowstuffinmargin
\showstuffinmarginfalse
\def\jjtags{\ifx\shlhetal\relax 
  \else
\ifx\shlhetal\undefinedcontrolseq
\else
\showstuffinmargintrue
\ifx\all\relax\else\expandafter\t@gsjj\all,\fi\fi \fi
}

\def\tagjj#1{\realstag{#1}\mginpar{\zeigen{#1}}}
\def\citejj#1{\rechnen{#1}\mginpar{\zeigen{#1}}}     

\def\rechnen#1{\expandafter\ifx\csname stag#1\endcsname\relax ??\else
                           \csname stag#1\endcsname\fi}

\newdimen\theight

\def\marginfont{\sevenrm}

\def\trymarginbox#1{\setbox0=\hbox{\marginfont\hskip\marginshift #1}%
		\global\marginshift\wd0 
		\global\advance\marginshift\margindelta}

\def \mginpar#1{%
\ifvmode\setbox0\hbox to \hsize{\hfill\rlap{\marginfont\quad#1}}%
\ht0 0cm
\dp0 0cm
\box0\vskip-\baselineskip
\else 
             \vadjust{\trymarginbox{#1}%
		\ifdim\marginshift>\marginmax \global\marginshift\marginmin
			\trymarginbox{#1}%
                \fi
             \theight=\ht0
             \advance\theight by \dp0    \advance\theight by \lineskip
             \kern -\theight \vbox to \theight{\rightline{\rlap{\box0}}%
\vss}}\fi}


\def\t@gsoff#1,{\def\@{#1}\ifx\@\empty\let\next=\relax\else\let\next=\t@gsoff
   \def\@@{p}\ifx\@\@@\else
   \expandafter\gdef\csname#1cite\endcsname##1{\zeigen{##1}}
   \expandafter\gdef\csname#1page\endcsname##1{?}
   \expandafter\gdef\csname#1tag\endcsname##1{\zeigen{##1}}\fi\fi\next}
\def\verbatimtags{\showstuffinmarginfalse
\ifx\all\relax\else\expandafter\t@gsoff\all,\fi}
\def\zeigen#1{\hbox{$\langle$}#1\hbox{$\rangle$}}


\def\margintag#1{\ifshowstuffinmargin\mginpar{\zeigen{#1}}\fi}

\def\(#1){\edef\dot@g{\ifmmode\ifinner(\hbox{\noexpand\etag{#1}})
   \else\noexpand\eqno(\hbox{\noexpand\etag{#1}})\fi
   \else(\noexpand\ecite{#1})\fi}\dot@g}

\newif\ifbr@ck
\def\eat#1{}
\def\[#1]{\br@cktrue[\br@cket#1'X]}
\def\br@cket#1'#2X{\def\temp{#2}\ifx\temp\empty\let\next\eat
   \else\let\next\br@cket\fi
   \ifbr@ck\br@ckfalse\br@ck@t#1,X\else\br@cktrue#1\fi\next#2X}
\def\br@ck@t#1,#2X{\def\temp{#2}\ifx\temp\empty\let\neext\eat
   \else\let\neext\br@ck@t\def\temp{,}\fi
   \def\teemp{#1}\ifx\teemp\empty\else\rcite{#1}\fi\temp\neext#2X}
\def\resetbr@cket{\gdef\[##1]{[\rtag{##1}]}}
\def\references{\resetbr@cket\newsection References\par}

\newtoks\symb@ls\newtoks\s@mb@ls\newtoks\p@gelist\n@wcount\ftn@mber
    \ftn@mber=1\newif\ifftn@mbers\ftn@mbersfalse\newif\ifbyp@ge\byp@gefalse
\def\defm@rk{\ifftn@mbers\n@mberm@rk\else\symb@lm@rk\fi}
\def\n@mberm@rk{\xdef\m@rk{{\the\ftn@mber}}%
    \global\advance\ftn@mber by 1 }
\def\rot@te#1{\let\temp=#1\global#1=\expandafter\r@t@te\the\temp,X}
\def\r@t@te#1,#2X{{#2#1}\xdef\m@rk{{#1}}}
\def\b@@st#1{{$^{#1}$}}\def\str@p#1{#1}
\def\symb@lm@rk{\ifbyp@ge\rot@te\p@gelist\ifnum\expandafter\str@p\m@rk=1 
    \s@mb@ls=\symb@ls\fi\write\f@nsout{\number\count0}\fi \rot@te\s@mb@ls}
\def\byp@ge{\byp@getrue\n@wwrite\f@nsin\openin\f@nsin=\jobname.fns 
    \n@wcount\currentp@ge\currentp@ge=0\p@gelist={0}
    \re@dfns\closein\f@nsin\rot@te\p@gelist
    \n@wread\f@nsout\openout\f@nsout=\jobname.fns }
\def\m@kelist#1X#2{{#1,#2}}
\def\re@dfns{\ifeof\f@nsin\let\next=\relax\else\read\f@nsin to \f@nline
    \ifx\f@nline\v@idline\else\let\t@mplist=\p@gelist
    \ifnum\currentp@ge=\f@nline
    \global\p@gelist=\expandafter\m@kelist\the\t@mplistX0
    \else\currentp@ge=\f@nline
    \global\p@gelist=\expandafter\m@kelist\the\t@mplistX1\fi\fi
    \let\next=\re@dfns\fi\next}
\def\symbols#1{\symb@ls={#1}\s@mb@ls=\symb@ls} 
\def\bigsymbol{\textstyle}
\symbols{\bigsymbol\ast,\dagger,\ddagger,\sharp,\flat,\natural,\star}
\def\ftnumbers{\ftn@mberstrue} \def\ftsymbols{\ftn@mbersfalse}
\def\paginal{\byp@ge} \def\resetftnumbers{\ftn@mber=1}
\def\ftnote#1{\defm@rk\expandafter\expandafter\expandafter\footnote
    \expandafter\b@@st\m@rk{#1}}

\long\def\jump#1\endjump{}
\def\ssum{\mathop{\lower .1em\hbox{$\textstyle\Sigma$}}\nolimits}

\def\qed{\nobreak\kern 1em \vrule height .5em width .5em depth 0em}
\def\newneq{\hbox{\rlap{\hbox to 1\wd9{\hss$=$\hss}}\raise .1em 
   \hbox to 1\wd9{\hss$\scriptscriptstyle/$\hss}}}
\def\subsetne{\setbox9 = \hbox{$\subset$}\mathrel{\hbox{\rlap
   {\lower .4em \newneq}\raise .13em \hbox{$\subset$}}}}
\def\supsetne{\setbox9 = \hbox{$\subset$}\mathrel{\hbox{\rlap
   {\lower .4em \newneq}\raise .13em \hbox{$\supset$}}}}

\def\vbar{\mathchoice{\vrule height6.3ptdepth-.5ptwidth.8pt\kern-.8pt}
   {\vrule height6.3ptdepth-.5ptwidth.8pt\kern-.8pt}
   {\vrule height4.1ptdepth-.35ptwidth.6pt\kern-.6pt}
   {\vrule height3.1ptdepth-.25ptwidth.5pt\kern-.5pt}}
\def\f@dge{\mathchoice{}{}{\mkern.5mu}{\mkern.8mu}}
\def\b@c#1#2{{\rm \mkern#2mu\vbar\mkern-#2mu#1}}
\def\b@b#1{{\rm I\mkern-3.5mu #1}}
\def\b@a#1#2{{\rm #1\mkern-#2mu\f@dge #1}}
\def\bb#1{{\count4=`#1 \advance\count4by-64 \ifcase\count4\or\b@a A{11.5}\or
   \b@b B\or\b@c C{5}\or\b@b D\or\b@b E\or\b@b F \or\b@c G{5}\or\b@b H\or
   \b@b I\or\b@c J{3}\or\b@b K\or\b@b L \or\b@b M\or\b@b N\or\b@c O{5} \or
   \b@b P\or\b@c Q{5}\or\b@b R\or\b@a S{8}\or\b@a T{10.5}\or\b@c U{5}\or
   \b@a V{12}\or\b@a W{16.5}\or\b@a X{11}\or\b@a Y{11.7}\or\b@a Z{7.5}\fi}}

\catcode`\X=11 \catcode`\@=12




\let\thischap\jobname

\def\partof#1{\csname returnthe#1part\endcsname}
\def\chapof#1{\csname returnthe#1chap\endcsname}

\def\setchapter#1,#2,#3.{%
  \expandafter\def\csname returnthe#1part\endcsname{#2}%
  \expandafter\def\csname returnthe#1chap\endcsname{#3}%
}

\setchapter 300a,A,I.
\setchapter 300b,A,II.
\setchapter 300c,A,III.
\setchapter 300d,A,IV.
\setchapter 300e,A,V.
\setchapter 300f,A,VI.
\setchapter 300g,A,VII.
\setchapter   88,B,I.
\setchapter  600,B,II.
\setchapter  705,B,III.

\def\cprefix#1{
\edef\theotherpart{\partof{#1}}\edef\theotherchap{\chapof{#1}}%
\ifx\theotherpart\thispart
   \ifx\theotherchap\thischap 
    \else 
     \theotherchap%
    \fi
   \else 
     \theotherpart.\theotherchap\fi}

\def\sectioncite[#1]#2{%
     \cprefix{#2}#1}

\edef\thispart{\partof{\thischap}}
\edef\thischap{\chapof{\thischap}}



\expandafter\ifx\csname citeadd.tex\endcsname\relax
\expandafter\gdef\csname citeadd.tex\endcsname{}
\else \message{Hey!  Apparently you were trying to
\string\input{citeadd.tex} twice.   This does not make sense.} 
\errmessage{Please edit your file (probably \jobname.tex) and remove
any duplicate ``\string\input'' lines}\endinput\fi

\def\sciteu{\sciteerror{undefined}}

\def\sciteerror#1#2{{\mathortextbf{\scite{#2}}}\complainaboutcitation{#1}{#2}}
\def\mathortextbf#1{\hbox{\bf #1}}
\def\complainaboutcitation#1#2{%
\vadjust{\line{\llap{---$\!\!>$ }\qquad scite$\{$#2$\}$ #1\hfil}}}

\sectno=-1   
\localtags
\verbatimtags
\NoBlackBoxes
\define\mr{\medskip\roster}
\define\sn{\smallskip\noindent}
\define\mn{\medskip\noindent}
\define\bn{\bigskip\noindent}
\define\ub{\underbar}

\define\ermn{\endroster\medskip\noindent}

\define\dbcu{\dsize\bigcup}
\define \nl{\newline}
\magnification=\magstep 1
\documentstyle {amsppt}

\expandafter\ifx\csname alice2jlem.tex\endcsname\relax
  \expandafter\xdef\csname alice2jlem.tex\endcsname{\the\catcode`@}
\else \message{Hey!  Apparently you were trying to
\string\input{alice2jlem.tex}  twice.   This does not make sense.}
\errmessage{Please edit your file (probably \jobname.tex) and remove
any duplicate ``\string\input'' lines}\endinput\fi

\expandafter\ifx\csname bib4plain.tex\endcsname\relax
  \expandafter\gdef\csname bib4plain.tex\endcsname{}
\else \message{Hey!  Apparently you were trying to \string\input
  bib4plain.tex twice.   This does not make sense.}
\errmessage{Please edit your file (probably \jobname.tex) and remove
any duplicate ``\string\input'' lines}\endinput\fi

\def\renewcommand{\newcommand}	       
\edef\cite{\the\catcode`@}%
\catcode`@ = 11
\let\@oldatcatcode = \cite
\chardef\@letter = 11
\chardef\@other = 12
%
%
%
%
\def\@innerdef#1#2{\edef#1{\expandafter\noexpand\csname #2\endcsname}}%
%
%
\@innerdef\@innernewcount{newcount}%
\@innerdef\@innernewdimen{newdimen}%
\@innerdef\@innernewif{newif}%
\@innerdef\@innernewwrite{newwrite}%
%
%
%
\def\@gobble#1{}%
%
%
%
\ifx\inputlineno\@undefined
   \let\@linenumber = \empty 
\else
   \def\@linenumber{\the\inputlineno:\space}%
\fi
%
%
%
\def\@futurenonspacelet#1{\def\cs{#1}%
   \afterassignment\@stepone\let\@nexttoken=
}%
\begingroup 
\def\\{\global\let\@stoken= }%
\\ 
\endgroup
\def\@stepone{\expandafter\futurelet\cs\@steptwo}%
\def\@steptwo{\expandafter\ifx\cs\@stoken\let\@@next=\@stepthree
   \else\let\@@next=\@nexttoken\fi \@@next}%
\def\@stepthree{\afterassignment\@stepone\let\@@next= }%
%
%
%
\def\@getoptionalarg#1{%
   \let\@optionaltemp = #1%
   \let\@optionalnext = \relax
   \@futurenonspacelet\@optionalnext\@bracketcheck
}%
%
%
\def\@bracketcheck{%
   \ifx [\@optionalnext
      \expandafter\@@getoptionalarg
   \else
      \let\@optionalarg = \empty
      \expandafter\@optionaltemp
   \fi
}%
\def\@@getoptionalarg[#1]{%
   \def\@optionalarg{#1}%
   \@optionaltemp
}%
%
%
%
\def\@nnil{\@nil}%
\def\@fornoop#1\@@#2#3{}%
\def\@for#1:=#2\do#3{%
   \edef\@fortmp{#2}%
   \ifx\@fortmp\empty \else
      \expandafter\@forloop#2,\@nil,\@nil\@@#1{#3}%
   \fi
}%
\def\@forloop#1,#2,#3\@@#4#5{\def#4{#1}\ifx #4\@nnil \else
       #5\def#4{#2}\ifx #4\@nnil \else#5\@iforloop #3\@@#4{#5}\fi\fi
}%
\def\@iforloop#1,#2\@@#3#4{\def#3{#1}\ifx #3\@nnil
       \let\@nextwhile=\@fornoop \else
      #4\relax\let\@nextwhile=\@iforloop\fi\@nextwhile#2\@@#3{#4}%
}%
%
%
%
\@innernewif\if@fileexists
\def\@testfileexistence{\@getoptionalarg\@finishtestfileexistence}%
\def\@finishtestfileexistence#1{%
   \begingroup
      \def\extension{#1}%
      \immediate\openin0 =
         \ifx\@optionalarg\empty\jobname\else\@optionalarg\fi
         \ifx\extension\empty \else .#1\fi
         \space
      \ifeof 0
         \global\@fileexistsfalse
      \else
         \global\@fileexiststrue
      \fi
      \immediate\closein0
   \endgroup
}%
%
%
%
%
\def\bibliographystyle#1{%
   \@readauxfile
   \@writeaux{\string\bibstyle{#1}}%
}%
\let\bibstyle = \@gobble
%
%
\let\bblfilebasename = \jobname
\def\bibliography#1{%
   \@readauxfile
   \@writeaux{\string\bibdata{#1}}%
   \@testfileexistence[\bblfilebasename]{bbl}%
   \if@fileexists
      \nobreak
      \@readbblfile
   \fi
}%
\let\bibdata = \@gobble
%
%
\def\nocite#1{%
   \@readauxfile
   \@writeaux{\string\citation{#1}}%
}%
\@innernewif\if@notfirstcitation
%
%
\def\cite{\@getoptionalarg\@cite}%
%
%
\def\@cite#1{%
   \let\@citenotetext = \@optionalarg
   \printcitestart
   \nocite{#1}%
   \@notfirstcitationfalse
   \@for \@citation :=#1\do
   {%
      \expandafter\@onecitation\@citation\@@
   }%
   \ifx\empty\@citenotetext\else
      \printcitenote{\@citenotetext}%
   \fi
   \printcitefinish
}%
\newif\ifweareinprivate
\weareinprivatetrue
\ifx\shlhetal\undefinedcontrolseq\weareinprivatefalse\fi
\ifx\shlhetal\relax\weareinprivatefalse\fi
\def\@onecitation#1\@@{%
   \if@notfirstcitation
      \printbetweencitations
   \fi
   \expandafter \ifx \csname\@citelabel{#1}\endcsname \relax
      \if@citewarning
         \message{\@linenumber Undefined citation `#1'.}%
      \fi
     \ifweareinprivate
      \expandafter\gdef\csname\@citelabel{#1}\endcsname{%
\strut 
\vadjust{\vskip-\dp\strutbox
\vbox to 0pt{\vss\parindent0cm \leftskip=\hsize 
\advance\leftskip3mm
\advance\hsize 4cm\strut\openup-4pt 
\rightskip 0cm plus 1cm minus 0.5cm ?  #1 ?\strut}}
         {\tt
            \escapechar = -1
            \nobreak\hskip0pt\pfeilsw
            \expandafter\string\csname#1\endcsname
             \pfeilso
            \nobreak\hskip0pt
         }%
      }%
     \else  
      \expandafter\gdef\csname\@citelabel{#1}\endcsname{%
            {\tt\expandafter\string\csname#1\endcsname}
      }%
     \fi  
   \fi
   \csname\@citelabel{#1}\endcsname
   \@notfirstcitationtrue
}%
%
%
\def\@citelabel#1{b@#1}%
%
%
\def\@citedef#1#2{\expandafter\gdef\csname\@citelabel{#1}\endcsname{#2}}%
%
%
%
\def\@readbblfile{%
   \ifx\@itemnum\@undefined
      \@innernewcount\@itemnum
   \fi
   \begingroup
      \def\begin##1##2{%
         \setbox0 = \hbox{\biblabelcontents{##2}}%
         \biblabelwidth = \wd0
      }%
      \def\end##1{}
      %
      %
      \@itemnum = 0
      \def\bibitem{\@getoptionalarg\@bibitem}%
      \def\@bibitem{%
         \ifx\@optionalarg\empty
            \expandafter\@numberedbibitem
         \else
            \expandafter\@alphabibitem
         \fi
      }%
      \def\@alphabibitem##1{%
         \expandafter \xdef\csname\@citelabel{##1}\endcsname {\@optionalarg}%
         \ifx\biblabelprecontents\@undefined
            \let\biblabelprecontents = \relax
         \fi
         \ifx\biblabelpostcontents\@undefined
            \let\biblabelpostcontents = \hss
         \fi
         \@finishbibitem{##1}%
      }%
      \def\@numberedbibitem##1{%
         \advance\@itemnum by 1
         \expandafter \xdef\csname\@citelabel{##1}\endcsname{\number\@itemnum}%
         \ifx\biblabelprecontents\@undefined
            \let\biblabelprecontents = \hss
         \fi
         \ifx\biblabelpostcontents\@undefined
            \let\biblabelpostcontents = \relax
         \fi
         \@finishbibitem{##1}%
      }%
      \def\@finishbibitem##1{%
         \biblabelprint{\csname\@citelabel{##1}\endcsname}%
         \@writeaux{\string\@citedef{##1}{\csname\@citelabel{##1}\endcsname}}%
         \ignorespaces
      }%
      %
      %
      \let\em = \bblem
      \let\newblock = \bblnewblock
      \let\sc = \bblsc
      \frenchspacing
      \clubpenalty = 4000 \widowpenalty = 4000
      \tolerance = 10000 \hfuzz = .5pt
      \everypar = {\hangindent = \biblabelwidth
                      \advance\hangindent by \biblabelextraspace}%
      \bblrm
      \parskip = 1.5ex plus .5ex minus .5ex
      \biblabelextraspace = .5em
      \bblhook
      \input \bblfilebasename.bbl
   \endgroup
}%
%
%
\@innernewdimen\biblabelwidth
\@innernewdimen\biblabelextraspace
%
%
%
\def\biblabelprint#1{%
   \noindent
   \hbox to \biblabelwidth{%
      \biblabelprecontents
      \biblabelcontents{#1}%
      \biblabelpostcontents
   }%
   \kern\biblabelextraspace
}%
%
%
%
\def\biblabelcontents#1{{\bblrm [#1]}}%
%
%
\def\bblrm{\rm}%
%
%
\def\bblem{\it}%
%
%
\def\bblsc{\ifx\@scfont\@undefined
              \font\@scfont = cmcsc10
           \fi
           \@scfont
}%
%
%
\def\bblnewblock{\hskip .11em plus .33em minus .07em }%
%
%
\let\bblhook = \empty
%
%
%
\def\printcitestart{[}
\def\printcitefinish{]}
\def\printbetweencitations{, }
\def\printcitenote#1{, #1}
%
%
%
\let\citation = \@gobble
%
%
%
\@innernewcount\@numparams
%
%
\def\newcommand#1{%
   \def\@commandname{#1}%
   \@getoptionalarg\@continuenewcommand
}%
%
%
\def\@continuenewcommand{%
   \@numparams = \ifx\@optionalarg\empty 0\else\@optionalarg \fi \relax
   \@newcommand
}%
%
%
\def\@newcommand#1{%
   \def\@startdef{\expandafter\edef\@commandname}%
   \ifnum\@numparams=0
      \let\@paramdef = \empty
   \else
      \ifnum\@numparams>9
         \errmessage{\the\@numparams\space is too many parameters}%
      \else
         \ifnum\@numparams<0
            \errmessage{\the\@numparams\space is too few parameters}%
         \else
            \edef\@paramdef{%
               \ifcase\@numparams
                  \empty  No arguments.
               \or ####1%
               \or ####1####2%
               \or ####1####2####3%
               \or ####1####2####3####4%
               \or ####1####2####3####4####5%
               \or ####1####2####3####4####5####6%
               \or ####1####2####3####4####5####6####7%
               \or ####1####2####3####4####5####6####7####8%
               \or ####1####2####3####4####5####6####7####8####9%
               \fi
            }%
         \fi
      \fi
   \fi
   \expandafter\@startdef\@paramdef{#1}%
}%
%
%
%
%
\def\@readauxfile{%
   \if@auxfiledone \else 
      \global\@auxfiledonetrue
      \@testfileexistence{aux}%
      \if@fileexists
         \begingroup
            \endlinechar = -1
            \catcode`@ = 11
            \input \jobname.aux
         \endgroup
      \else
         \message{\@undefinedmessage}%
         \global\@citewarningfalse
      \fi
      \immediate\openout\@auxfile = \jobname.aux
   \fi
}%
%
%
\newif\if@auxfiledone
\ifx\noauxfile\@undefined \else \@auxfiledonetrue\fi
%
%
%
%
\@innernewwrite\@auxfile
\def\@writeaux#1{\ifx\noauxfile\@undefined \write\@auxfile{#1}\fi}%
%
%
%
\ifx\@undefinedmessage\@undefined
   \def\@undefinedmessage{No .aux file; I won't give you warnings about
                          undefined citations.}%
\fi
%
%
\@innernewif\if@citewarning
\ifx\noauxfile\@undefined \@citewarningtrue\fi
%
%
%
\catcode`@ = \@oldatcatcode

\def\pfeilso{\leavevmode
            \vrule width 1pt height9pt depth 0pt\relax
           \vrule width 1pt height8.7pt depth 0pt\relax
           \vrule width 1pt height8.3pt depth 0pt\relax
           \vrule width 1pt height8.0pt depth 0pt\relax
           \vrule width 1pt height7.7pt depth 0pt\relax
            \vrule width 1pt height7.3pt depth 0pt\relax
            \vrule width 1pt height7.0pt depth 0pt\relax
            \vrule width 1pt height6.7pt depth 0pt\relax
            \vrule width 1pt height6.3pt depth 0pt\relax
            \vrule width 1pt height6.0pt depth 0pt\relax
            \vrule width 1pt height5.7pt depth 0pt\relax
            \vrule width 1pt height5.3pt depth 0pt\relax
            \vrule width 1pt height5.0pt depth 0pt\relax
            \vrule width 1pt height4.7pt depth 0pt\relax
            \vrule width 1pt height4.3pt depth 0pt\relax
            \vrule width 1pt height4.0pt depth 0pt\relax
            \vrule width 1pt height3.7pt depth 0pt\relax
            \vrule width 1pt height3.3pt depth 0pt\relax
            \vrule width 1pt height3.0pt depth 0pt\relax
            \vrule width 1pt height2.7pt depth 0pt\relax
            \vrule width 1pt height2.3pt depth 0pt\relax
            \vrule width 1pt height2.0pt depth 0pt\relax
            \vrule width 1pt height1.7pt depth 0pt\relax
            \vrule width 1pt height1.3pt depth 0pt\relax
            \vrule width 1pt height1.0pt depth 0pt\relax
            \vrule width 1pt height0.7pt depth 0pt\relax
            \vrule width 1pt height0.3pt depth 0pt\relax}

\def\pfeilsw{ \leavevmode 
            \vrule width 1pt height0.3pt depth 0pt\relax
            \vrule width 1pt height0.7pt depth 0pt\relax
            \vrule width 1pt height1.0pt depth 0pt\relax
            \vrule width 1pt height1.3pt depth 0pt\relax
            \vrule width 1pt height1.7pt depth 0pt\relax
            \vrule width 1pt height2.0pt depth 0pt\relax
            \vrule width 1pt height2.3pt depth 0pt\relax
            \vrule width 1pt height2.7pt depth 0pt\relax
            \vrule width 1pt height3.0pt depth 0pt\relax
            \vrule width 1pt height3.3pt depth 0pt\relax
            \vrule width 1pt height3.7pt depth 0pt\relax
            \vrule width 1pt height4.0pt depth 0pt\relax
            \vrule width 1pt height4.3pt depth 0pt\relax
            \vrule width 1pt height4.7pt depth 0pt\relax
            \vrule width 1pt height5.0pt depth 0pt\relax
            \vrule width 1pt height5.3pt depth 0pt\relax
            \vrule width 1pt height5.7pt depth 0pt\relax
            \vrule width 1pt height6.0pt depth 0pt\relax
            \vrule width 1pt height6.3pt depth 0pt\relax
            \vrule width 1pt height6.7pt depth 0pt\relax
            \vrule width 1pt height7.0pt depth 0pt\relax
            \vrule width 1pt height7.3pt depth 0pt\relax
            \vrule width 1pt height7.7pt depth 0pt\relax
            \vrule width 1pt height8.0pt depth 0pt\relax
            \vrule width 1pt height8.3pt depth 0pt\relax
            \vrule width 1pt height8.7pt depth 0pt\relax
            \vrule width 1pt height9pt depth 0pt\relax
      }


\def\widestnumber#1#2{}

\def\citewarning#1{\ifx\shlhetal\relax 
    \else
    \par{#1}\par
    \fi
}

\def\rm{\fam0 \tenrm}

\def\fakesubhead#1\endsubhead{\bigskip\noindent{\bf#1}\par}



%
%
%

%

\font\textrsfs=rsfs10
\font\scriptrsfs=rsfs7
\font\scriptscriptrsfs=rsfs5

\newfam\rsfsfam
\textfont\rsfsfam=\textrsfs
\scriptfont\rsfsfam=\scriptrsfs
\scriptscriptfont\rsfsfam=\scriptscriptrsfs

\edef\oldcatcodeofat{\the\catcode`\@}
\catcode`\@11

\def\Cal@@#1{\noaccents@ \fam \rsfsfam #1}

\catcode`\@\oldcatcodeofat


\expandafter\ifx \csname margininit\endcsname \relax\else\margininit\fi

\long\def\red#1\endred{}
\long\def\green#1\endgreen{}
\long\def\blue#1\endblue{}

\def\endred{ \unmatched endred! }
\def\endgreen{ \unmatched endgreen! }
\def\endblue{ \unmatched endblue! }

\ifx\latexcolors\undefinedcs\def\latexcolors{}\fi

\def\emptycs{}
\def\evaluatelatexcolors{%
        \ifx\latexcolors\emptycs\else
        \expandafter\xxevaluate\latexcolors\xxfertig\evaluatelatexcolors\fi}
\def\xxevaluate#1,#2\xxfertig{\setupthiscolor{#1}%
        \def\latexcolors{#2}}

\font\smallfont=cmsl7
\def\rutgerscolor{\ifmmode\else\endgraf\fi\smallfont
\advance\leftskip0.5cm\relax}
\def\setupthiscolor#1{\edef\tmptmpcs{\noexpand\bgroup\noexpand\rutgerscolor
\noexpand\def\noexpand\currentcolor{#1}%
\noexpand}%
\expandafter\let\csname#1\endcsname\tmptmpcs
\def\tmptmpcs{\checkColorUnmatched{#1}\popthecolor}
\expandafter\let\csname end#1\endcsname\tmptmpcs}

\def\checkColorUnmatched#1{\def\expectcolor{#1}%
    \ifx\expectcolor\currentcolor   
    \else \edef\failhere{\noexpand\tryingToClose '\currentcolor' with end\expectcolor}\failhere\fi}

\def\currentcolor{???}

\def\popthecolor{\ifmmode\else\endgraf\fi\egroup}

\expandafter\def\csname#1\endcsname{}

\evaluatelatexcolors

 \let\outerhead\head
 \def\head{\innerhead}
 \let\innerhead\outerhead

 \let\outersubhead\subhead
 \def\subhead{\innersubhead}
 \let\innersubhead\outersubhead

 \let\outersubsubhead\subsubhead
 \def\subsubhead{\innersubsubhead}
 \let\innersubsubhead\outersubsubhead

 \def\proclaim{\innerproclaim}
 \let\innerproclaim\outerproclaim

 %
 %
 %
 %

\def\demo#1{\medskip\noindent{\it #1.\/}}
\def\enddemo{\smallskip}

\def\remark#1{\medskip\noindent{\it #1.\/}}
\def\endremark{\smallskip}

\topmatter
\title{FPL may be equivalent to FO but not equivalent to PFP} \endtitle
\author {Saharon Shelah \thanks {\null\newline I would like to thank 
Alice Leonhardt for the beautiful typing. Publication E31.\null\newline
} \endthanks} \endauthor 

\affil {Institute of Mathematics\\
 The Hebrew University\\
 Jerusalem, Israel
 \medskip
 Rutgers University\\
 Mathematics Department\\
 New Brunswick, NJ  USA} \endaffil
\bn
\abstract  Given a class of finite models we would like to expand each
model (allowing new elements but the old universe is a separate sort), making
the expressive power of LFP (least fix point logic) and PFP (inductive logic)
similar while not changing the expressive power of FO (first order logic).
This continues in \cite{GISh:525}. \endabstract
\endtopmatter
\document  
\newpage

\head {Anotated Content} \endhead  \resetall \resetall
\bn
\S1 \ub{The Construction}
\mr
\item "{{}}"  [We deal with the construction used, basically for coding
$\bar b \in S_\ell[M_i]$ in $M^2_i$ we use an $E_1$-equivalent class, in it
for each $a \in M_i$ we attach a set $B_{\langle \ell,\bar b,a \rangle}$,
size determined by $i$ and $\ell$, and on it put a graph coding the answer
if $a \notin \text{ Rang}(\bar b)$ and coding $b_j$ if $a = b_j$ (we can
assume $\|M_i\| > 2 \ell g(\bar b) > 100$).  This is simpler than the previous
random graph.  Still inside we have the cliques of size $h(e)$ and the 
cycles.
We show that F.O. says nothing new on $M_i$, while 
LFP logic expresses all the
global relations attached.  Also we essentially eliminate quantifiers showing
what $PFP$ in $M^2_i$ can say when restricted to $M_i$, the original model.
So we can get: $LFP, PFP$ expressing the same in $\{M^2_i:i \in \Bbb N\}$ 
when restricted to $\{M_i:i \in \Bbb N\}$.  
But it is not shown: $LFP, PFP$ are
equivalent on $\{M^2_i:i\}$ (i.e. closure under isomorphism).]
\endroster
\bn
\S2 \ub{Discussion}
\mr
\item "{{}}"  [Discuss some points.]
\endroster
\bn
\S3 \ub{Redoing}
\mr
\item "{{}}"  [Here we redo \S1 somewhat differently to get the desired
original question.  For this $2^{h(i,\ell)} < h(i,\ell + 1)$ is really used.
The main change is that the $S_\ell(M_i)$ is not a global $m(\ell)$-place
relation on $M_i$, \ub{but} a global $m(\ell)$-place function from $M_i$ to
$\Bbb N \restriction h_2(i,\ell +3)$.  Probably \ub{more natural} is to ask
$S_i(M_i)$ is a global $m(i)$ place function from $M_i$ to
$\Bbb N \restriction h_3(M,\ell)$ where 
$h_3(M,\ell) \le h_2(i,\ell)$ but does not matter.  We indicate the changes.]
\endroster
\newpage

\head{\S1 The Construction} \endhead  \resetall \resetall
\bigskip

\demo{\stag{1.1} Hypothesis}
\mr
\item "{$(a)$}"  $\tau$ a finite vocabulary
\sn
\item "{$(b)$}"  $K$ a class of finite $\tau$-models, 
\sn
\item "{$(c)$}"  $\langle S_i:i < \omega \rangle$ a list of global relations
on $K$, or even on $K_1 = \{M:$ for some $i < \omega$ we have 
$M,M^1_i\}$ (see below), $S_i$ has $m(S_i)$-places, without loss of
generality $0 < m(S_i)$, for simplicity $S_i$ is irreflexive and lastly
$\langle M_j:j < \omega \rangle$ list $K$ with no two isomorphic, i.e. 
this list a set of representatives
\sn
\item "{$(d)$}"  $h:\omega \rightarrow \omega$ goes to infinite (slowly).
\endroster
\enddemo
\bn
\margintag{1.2}\ub{\stag{1.2} Choice}:  For $k,t,m < \omega$ let $G = G_{k,t,m}$ ($k \ge 10$
for simplicity so that $Y$ below exists) be a graph such that we have
$x_{\ell,j} \in G$ (for $\ell < k,j < m$, with no repetitions) and
\mr
\item "{$(a)$}"  $\{x_{\ell,j}:\ell < k\}$ is a maximal clique
\sn
\item "{$(b)$}"  if $j_1 = j_2 + 1$ mod $m$ then $\langle x_{\ell,j_1}:\ell <
k \rangle \char 94 \langle x_{\ell,j_2}:\ell < k \rangle$ has a special
property (say for $\ell_1,\ell_2 < k:[\{x_{\ell_1,j_1},x_{\ell_2,j}\}$ an edge
$\Leftrightarrow (\ell_1,\ell_1 + k) \in Y]$ where $Y \subseteq \{0,\dotsc,
k-1\} \times \{k,k+1,\dotsc,2k-1\}$) is random enough: 
$\forall x < 2k \exists yz
[(x,y) \in Y \wedge (x,z) \notin Y \wedge (z < k \equiv x \ge k)]$ and the
bipartite graph $(\{0,\dotsc,2k-1\},Y)$ is rigid
\sn
\item "{$(c)$}"  if $i_1 \ne i_2,i_2+1$ mod $m$, then $\langle
x_{\ell,i_1}:\ell < k \rangle \char 94 \langle x_{\ell,i_2}:\ell < k \rangle$
strongly fail the special property, say no 
$(x_{\ell_1,i_1},x_{\ell_2,i_2})$ is an edge
\sn
\item "{$(d)$}"  if $A,B \subseteq G$ are disjoint, $|A| + |B| < t$ and
$[i < m \Rightarrow A$ does not contain $\{x_{\ell,i}:\ell < k\}]$, then for
some $x \in G \backslash \{x_{\ell,j}:i,j\}$ we have 
$y \in A \Rightarrow yRx,y \in B \Rightarrow
\neg yRx$.
\endroster
\bn
\margintag{1.3}\ub{\stag{1.3} Choice}:  We choose functions such that $h_0(i) < h(i),
h_0(i)$ goes to infinity, $h_1:\omega \rightarrow \omega,h_2(i,\ell)$ is
defined for $\ell < h_0(i)$, the values are $<$ log$_2$ log$_2 h(i)$, the
functions $i \mapsto h_2(i,0)$ goes to infinity, $10 \le h_1(\ell) < 
h_1(\ell +1) < h_2(i,\ell) < \text{ log}_2[h_2(i,\ell + 1)]$.
\bn
\margintag{1.4}\ub{\stag{1.4} Construction}:  We assume below no incidental equalities
between elements occurs.

For $i < \omega,M_i \in K$ is given see \scite{1.1}(c), let
\mr
\item "{$(a)$}"  $M^0_i = M_i$ and
\sn
\item "{$(b)$}"  $M^1_i = M^0_i + \Bbb N \restriction h(i)$ so e.g. we
have predicates $P_0,P^{M^1_i}_0 = M^0_i$ and $P_1,P^{M^1_i}_1 = 
(\Bbb N \restriction i)
\cup M^0_i$, for $R$ a predicate of $\tau$ let $R^{M^1_i} = R^{M^0_i}$ 
similarly for function symbols of $\tau$ and $0$ individual constant, 
$+ \restriction h(i),\times \restriction h(i)$ and $(x+1) \restriction h(i)$ 
are, of course, $0,+,\times$ and $x+1$ restricted to $h(i)$ and 
$P_0,P_1 \notin \tau$
\sn
\item "{$(c)$}"  $M^2_i$ is defined as follows: \nl
\ub{set of elements}:  $|M^1_i| \cup \dbcu_{\ell < h_0(i)} C_{i,\ell}$ where
$$
\align
B_{i,\ell} = \{\langle \ell,\bar b,a \rangle:\,&\bar b \in {}^{m(S_\ell)}
(|M^0_i|) \text{ with no repetition}, \\
  &\,\text{i.e. } \bar b \text{ a sequence of members of}\\
  &\,M^0_i \text{ of length } m(S_\ell),s_1 < s_2 \Rightarrow b_{s_1} \ne
b_{s_2} \\
  &\text{and } a \in M^0_i\}
\endalign
$$
for $u = \langle \ell,\bar b,a \rangle \in B_{i,\ell}$ we let $(k(u),t(u),
m(u))$ be
$$
k(u) = h_1(\ell)
$$

$$
m(u) = \cases m(S_i)h_2(i,\ell) &\text{ \ub{if} } \bar b \in 
S_\ell[M_i] \\
m(S_\ell)h_2(i,\ell) + j + 2 &\text{ \ub{if} } \bar b \notin S_\ell[M_1],
a = b_j (j \text{ unique}) \\
  m(S_\ell)h_2(i,\ell) + 1 &\text{ \ub{if} } \bar b \notin S_\ell[M_1],
a \notin \text{ Rang}(\bar b) \endcases
$$

$$
t(u) = 2k(u)m(u)
$$
and let $u = \langle \ell(u),\bar b(u),a(u) \rangle$ \nl
and lastly

$$
C_{i,\ell} = \{ \langle u,c \rangle:c \in G_{k(u),t(u),m(u)} \text{ and }
u \in B_{i,\ell}\}.
$$
\endroster
\bn
\ub{Relations}:  
\mr
\item "{$(\alpha)$}"  those of $M^1_i$ (interpret same way)
\sn
\item "{$(\beta)$}"  $F$ a partial one-place function, Dom$(F) = 
M^2_i \backslash M^1_i,f((u,c)) = a(u)(\in M^1_i)$
\sn
\item "{$(\gamma)$}"  $E_1$ is an equivalence relation on $M^2_i$ such
that $(u_1,c_1)E_1(u_2,c_2)$ iff $\ell(u_1) = \ell(u_2)$
\sn
\item "{$(\delta)$}"  $E_2$ is an equivalence relation on $M^2_i$ 
refining $E_1:(u_1,c_1)E_2(u_2,c_2)$ iff $\ell(u_1) = \ell(u_2) \and 
\bar b(u_1) = \bar b(u_2)$
\sn
\item "{$(\varepsilon)$}"  $E_3$ is the equivalence relation 
$(u_1,c_1),E_3(u_2,c_2)$ iff $u_1 = u_2$, so it refines $E_2$.
\ermn
(d) Let $M^3_i = M^2_i \restriction (M^0_i \cup \text{ Dom}(E^{M^2_i}))$ 
(i.e. we omit the arithmetic). \nl
(e) Let for $t \le h_0(i)$

$$
M^4_{i,t} = M^2_i \restriction (M^0_i \cup \text{ Dom}(E^{M^2_i}_0) \cup
\Bbb N \restriction \tau)
$$
\mn
So $M^4_{i,0} = M^3_i,M^4_{i,h_0(i)} = M^2$.
\bigskip

\definition{\stag{1.5} Definition}  $K_\ell = \{M:M \cong M^\ell_i
\text{ for some } i\}$ for $\ell = 0,1,2,3$, so $K_0 = K$.

$$
K_4 = \{M:M \cong_{M^2_i} M^\ell_{i,t} \text{ for some } t \le h_0(i)
\text{ and } i < u\}.
$$
\enddefinition
\bigskip

\proclaim{\stag{1.6} Claim}  1) For $\ell_0 < \ell_1 \le 2$ and first order
formula $\varphi(\bar x) \in {\Cal L}_{\tau(K_{\ell_1})}$ there is 
a first order
$\psi(\bar x) \in {\Cal L}_{\tau(K_{\ell_0})}$ such that $M \in K_{\ell_1}
\Rightarrow \varphi(\bar x)^M \restriction P_{\ell_0} = 
\psi(\bar x)^{(M \restriction P_{\ell_0} \restriction \tau(K_{\ell_0}))}$
equivalently $i < \omega \Rightarrow \varphi(\bar x)^{M^{\ell_1}_i} 
\restriction M^{\ell_0}_i = \psi(\bar x)^{M^{\ell_0}_i}$. \nl
2)  Moreover, for every first order formula 
$\varphi(\bar x,\bar y) \in {\Cal L}_{\tau(K_{\ell_0})}$ there is a first
order formula $\psi(\bar x,\bar z) \in {\Cal L}_{\tau(K_{\ell_0})}$
such that: if $M \in K_{\ell_1},\bar c \in {}^{\ell g(\bar x)}M$, then for
some $\bar d \in {}^{\ell g(\bar z)}M$ we have 
$\varphi(x,\bar c)^M \restriction
P_{\ell_0} = \psi(\bar x,\bar d)^{M \restriction P_{\ell_0} \restriction
\tau(K_{\ell_0})}$ equivalently if $i < \omega,\bar c \in
{}^{\ell g(\bar y)}(M^{\ell_1}_i)$ then for some $\bar d \in
{}^{\ell g(\bar z)}(M^{\ell_0}_i)$

$$
\varphi(\bar x,\bar c)^{M^{\ell_1}_i} \restriction
M^{\ell_0}_i = \psi(\bar x,\bar d)^{M^{\ell_0}_i}.
$$
\mn
3)  We can replace $K_2$ by $K_4$.
\endproclaim
\bigskip

\demo{Proof}  Enough to prove for ($\ell_0,\ell_1) = \{(0,1),(1,2)\}$.
\mn
\ub{$(\ell_0,\ell_1) = (0,1)$}

By the addition theorem.
\mn
\ub{$(\ell_0,\ell_1) = (1,2)$}

Let $s$ be the quantifier depth of $\varphi$.  We first look at $M'_i = 
M^2_i \restriction \{x:x \in M^1_i \text{ or } x = (u,c),k(c) > s\}$, then
all the graphs we add to $M^1_i$ (to get $M'_i$) are $s$-random (as
$G_{k,t,m}$ is Min$\{k-1,m\}$-random).  Now find $\psi$ directly or e.g. 
we can first replace all the graphs by copying one such graph, then add 
isomorphisms between them, then replace it by one graph and a representative 
of each $E_2$-equivalence class, then inspect.

The rest: going to full $M^2_i$.  We can assume $h_2(i,0) > s$, noting f.o.
says little on circles.

Of course, we use: it is enough to find an equivalence formula for $i$
large enough. \nl
${{}}$ \hfill$\square_{\scite{1.6}}$
\enddemo
\bigskip

\proclaim{\stag{1.7} Claim}  1) For every $\ell$, for some $\varphi(\bar x)
\in LFP(\tau)$ we have $\varphi(\bar x)^{M^2_i} \restriction P_0 = 
\varphi(\bar x)^{M^2_i} \restriction M_i = S_\ell(M_i)$. \nl
2)  Similarly using $M^4_{i,\tau}$.
\endproclaim
\bigskip

\demo{Proof}  Let $\varphi(\bar x)$ say:

there is $y \in \text{ Dom}(E_1)$ such that:
\mr
\item "{$(a)$}"  in the graph $(y/E_3,R \restriction (y/E_3))$ there is
a maximal clique of size $h_1(\ell)$
\sn
\item "{$(b)$}"  there are $y_0,y_1,\dotsc,y_{m(S_i)-1}$ in $y/E_2$ pairwise
not $E_3$-equivalent (remember ``$S_i$ is irreflexive") and for
$j = 0,\dotsc,m(S_i)-1$, in the graph $(x_j/E_3,R \restriction (x_j/E_3))$ 
the set of maximal
cliques form a cycle of length $=$ mod $m(S_i)$ and $F(y_0) = x_0 \and
F(y_1) = x_1 \and \ldots F(y_{\ell g(\bar x-1)}) = x_{\ell g(\bar x-1)}$.
\ermn
This is easily expressible (well as alway we ignore the finitely many
$\{M_i:h_1(\ell) < h_0(i)\}$ as it can be corrected).
\hfill$\square_{\scite{1.7}}$
\enddemo
\bigskip

\proclaim{\stag{1.8} Claim}  1) If $\varphi(\bar x) \in PFP(\tau(K_2))$,
\ub{then} for some $\ell < \omega$ and \nl
$\psi(\bar x) \in PFP(\tau(K_1) \cup \{S_0,\dotsc,S_{\ell -1}\})$ we have \nl
$M_2 \in K_2 \Rightarrow \varphi(\bar x)^{M_2}
\restriction (P^{M_2}_1) = \psi(\bar x)^{M'}$ where \nl
$M_2 = (M_2 \restriction P_1 \restriction
\tau(K_1),S^{M_2 \restriction P_0 \restriction \tau(K_0)}_0,\dotsc,
S^{M_2 \restriction P_0 \restriction (K_0)}_{\ell -1})$. \nl
2) For any $\varphi(\bar x,\bar y) \in PFP(\tau(K_2))$ for some 
$\varphi(\bar x,\bar z) \in LFP(\tau(K_1))$ we have: if \nl
$M_2 \in K_2,\bar c \in {}^{\ell g(\bar x)}(M_2)$, \ub{then} for some
$\bar d \in {}^{\ell g(\bar z)}(M_2)$ we have \nl
$\varphi(\bar x,\bar c)^{M_2} \restriction P^{M_2}_1 = 
\psi(\bar x,\bar d)^{M_2 \restriction P_1 
\restriction \tau(K_1)} \restriction P$.
\endproclaim
\bigskip

\remark{\stag{1.8A} Remark}  If stuck in the middle, read \scite{3.8} 
which explicate more $R_i$'s role.
\endremark
\bigskip

\demo{Proof}  1), 2)  We can consider a number $n$, a $n$-place 
predicate-variable $Q$ and $\varphi = \varphi(Q,x_0,\dotsc,x_{n-1})$ 
first order in
the vocabulary of $\tau(K_2)$ and we would like to analyze for $M \in K_2,
\langle Q^M_{\varphi,v}:v < v^M_\varphi \rangle$ where $Q^M_{\varphi,0} =
\emptyset,Q^M_{\varphi,v+1} = \{\bar x:M_2 \models \varphi(Q^M_{\varphi,v},
\bar x)\}$; this will implicitly tell us how to choose $\varphi$ (for (1)
and for (2)).  let $k(*)$ be: (the quantifier depth of $\varphi) +n+2$.  
We look at $M \in K_2$.  So let $M = M^2_i$.  Let

$$
A^3_i = \{(u,c) \in M^2_i \backslash M^1_i:h_1(\ell(u)) \le k(*)\}
$$

$$
A^4_i = \{(u,c) \in M^2_i \backslash M^1_i:h_1(\ell(u)) \ge k(*)\}
$$
\mn
if $d \in A^3_i$ we know $(d/E_2,R \restriction (x/E_2))$ is isomorphic to
some $G$ as in \scite{1.2} and call the cycle $\langle x^d_{j,\ell}:j <
m^d,\ell < k^d \rangle$ so we demand $d'E_2d'' \Rightarrow \dsize
\bigwedge_{j,\ell} x^{d'}_{j,\ell} = x^{d''}_{j,\ell} \and (m^{d'},k^{d'}) =
(m^{d''},k^{d''})$.
\sn
Let us define a family ${\Cal R}_i = {\Cal R}(M^2_i)$ of $n$-place relations
$Q$ on $M^2_i$ such that: every function $f \in {\Cal F}_i = {\Cal F}
[M^2_i]$ (defined below) preserve $Q$, i.e.

$$
x_0,\dotsc,x_{n-1} \in \text{ Dom}(f) \Rightarrow Q(x_0,\dotsc,x_{n-1})
\equiv Q(f(x_0) \ldots f(x_{n-1}) 
$$
\mn
$\text{where}$
$$
\align
{\Cal F}_i = \biggl\{ f:&(a) \quad f \text{ is a partial permutation of }
M^2_i, \\
  &(b) \quad f \restriction M^1_i \text{ is the identity} \\
  &(c) \quad |\text{Dom }f| \le k(*) \\
  &(d) \quad \text{if } x \in (\text{Dom }f) \cap A^3_i \text{ then }
x E^{M^2_i}_0f(x) \\
  &(e) \quad f \text{ is a partial isomorphism (i.e. preserve } F(-),
\pm E_0,\pm E_1,\pm E_2) \\
  &(f) \quad \text{we can find } \langle j^*_x:x \in (\text{Dom }f) \cap
A^3_i \rangle \text{ such that } xE_2 x' \Leftrightarrow j^*_x = j^*_{x'}
\biggr\} \\
  &\text{ and if } x \in (\text{Dom }f) \cap A^3_i \text{ let } 
x_0 = x,x_1 = f(x), \text{ \ub{then} for} \\
  &s \in \{0,1\} \text{ in } x_s/E_3 
\text{ letting the cycle (see \scite{1.2}) consist of} \\ 
  &\langle x^s_{j,\ell}:\ell < k^s,j < m^s \rangle \text{ as } x_0e^{M^2_i}_0
f(x) \text{ necessarily} \\
  &m^0 = m^1,k^0 = k^1 (\le k(*)) \text{ now the demand is} \\
  &(\forall \ell < k^s)(\forall j < m^s)
[x_0 = x^0_{j,\ell} \leftrightarrow x_1 = s^1_{j+j_x,\ell}]
(\forall \ell < k^s)(\forall j < m^s) \\
  &[x_0Rx^0_{j,\ell} \leftrightarrow x_1Rx^1_{j+j,\ell}].
\endalign
$$
\mn
Remember $R$ is symmetric.
\enddemo
\bn
\margintag{1.8B}\ub{\stag{1.8B} Fact}:  $Q^{M^2_i}_{\varphi,v} \in {\Cal R}_i$.
\bigskip

\demo{Proof}  By induction on $v$.

Now each $Q \in {\Cal R}_i$ we can code by a relation (with more than $n$ places) on $M^1_i$, we need the $M'_i \restriction (P_j \backslash P_0)
= \Bbb N \restriction h_0(i)$ to code the $j^*_x$'s.
\enddemo
\bigskip

\definition{\stag{1.8C} Definition}   A coding of 
$Q \in {\Cal R}_i$ is giving the following information:
\mr
\item "{$(*)_1$}"  for $n_1 \le n$ and quantifier free complete type $r$
in the variable $x_{n_1},\dotsc,x_{n-1}$ and the predicate $E_0,E_1,E_2,R$,
the $n_1$-place relation
\sn
\item "{$(*)_2$}"  $\{\bar b \char 94 \bar c \in {}^{n_1}(M^1_i):
\ell g(\bar b) = n_1,\ell g(\bar c) = n-n_1 \text{ for some } \bar c' \in
{}^{n-n_1}(A^4_i)$, \nl
$\bar c \text{ realizes } \bar r, \bar b \char 94 \bar c \in
Q \text{ and } c_\ell = F(c'_\ell)\}$ (supercedes in margin)
\sn
\item "{$(*)_3$}"  for every division $\bar w = \langle w_0,w_1,w_2 \rangle$
of $\{0,\dotsc,n-1\}$ and quantifier free complete type $r=r(x_0,\dotsc,
x_{n-1})$ is $\tau(K_2)$ such that
$$
n' \in w_0 \Rightarrow P_1(x) \in r
$$

$$
n' \in w_1 \cup w_2 \Rightarrow \neg P_2(x) \in r
$$
\ermn

$$
\align
R^*_{\bar w,\bar r}(\bar Q) = \biggl\{ \bar b \char 94 \bar c:&\ell g(\bar b)
= n \text{ and for some } \bar b' \text{ we have} \\
  &(\alpha) \quad \ell g(\bar b') = n \\
  &(\beta) \quad m \in w_0 \Rightarrow b'_m = b_m \\
  &(\gamma) \quad m \in w_1 \Rightarrow b'_m = F(b_m) \and b_m \in A^3_i \\
  &(\delta) \quad m \in w_2 \Rightarrow b'_m = F(b_m) \and b_m \in A^4_i \\
  &(\varepsilon) \quad \bar c = \langle g_m:m \in w_1 \rangle \\
  &(\zeta) \quad \text{ if } j \in w_1, \text{ then } c_j \in P^{M^2_i}_1
\backslash P^{M^2_i}_0 (= M^1_i \backslash M^0_i) \\
  &\qquad \text{ code the following two sequences 
(of length } k^{b_m} \times m^{b_m}) \\
  &\qquad \langle \text{truth value } (b_j = x^{b_j}_{k,m}):k < k^{b_m},m <
m^{b_m} \rangle \\
  &\qquad \langle \text{truth value } (b_j Rx^{b_j}_{k,m}):k < k^{b_m},m <
m^{b_m} \rangle \biggr\}.
\endalign
$$
\enddefinition
\newpage

\head {\S2 Discussion} \endhead  \resetall \resetall
\bn
\margintag{2.1}\ub{\stag{2.1} Question}:  What is the size of $M^2_i$?
\bn
\ub{Answer}:  If answer is $< \|M_i\|^{h_0(i)}$, now $h_0(i)$ could grow
very slowly with $i$, so we are just barely above polynomial growth.  But
using Baire categoricity theorem and (finite) cardinality considerations, we
cannot do better.  Even if we like to code one global $n$-place relation 
on $K$, we may need to increase $M_i$ to $\sim \|M_i\|^{n/m}$ if $m$ is 
maximal arity of the new vocabulary.
\bn
\margintag{2.2}\ub{\stag{2.2} Question}:  Can we use $h_2(i,\ell) = h_2(i,0)$?
\bn
\ub{Answer}:  I find in \S1 no use of $h_2(i,\ell) < h_2(i,\ell +1)$ (or
$h_2(i,\ell) \le \text{ log}_2 h_2(i,\ell +1)$) though it had seemed to me
natural, however, in \S3 we shall use ?
\bn
\margintag{2.3}\ub{\stag{2.3} Question}:  Can we not add $\Bbb N \restriction h_0(i)$?
\bn
\ub{Answer}:  If some $S_i$ code a linear length $\ge h_0(i)$, \ub{then} yes.
\nl
Also if some $\varphi(Q,\bar x)$ is $a$ in the proof of \scite{1.8} and
$\varphi(Q,\bar x)$ is positive in $Q$,
the induction on $M_i$, takes $\ge h_0(i)$ steps, we can easily do this.
\newpage

\head {\S3 Redoing} \endhead  \resetall \resetall
\bigskip

We redo \S1, so 3.x is a variant of 1.x and we say the changes \nl
\ub{\sciteu{3.1}}  We replace clause (c) by:
\mr
\item "{$(c)$}"   $S_\ell$ is a global functional, i.e. gives for every
$M \in K$ a function $S^M_\ell$ from ${}^{m(S_\ell)}(M)$ to $\{0,\dotsc,
h_2(M,\ell -1)\},0 < h_2(M<\ell -1) < \omega$ and we have $h_2(M_i,\ell-3) =
h_2(i,\ell-1)$ (see later) and we let e.g. $h_2(M,-1) = 1 = h_2(M,-2) =
h_2(M,-3)$ and
{\roster
\itemitem{ $(*)$ }  if $f$ is an isomorphism from $M' \in K$ onto $M'' \in K$
\ub{then} \nl
$(S_\ell(M'))(x_0,\dotsc,x_{m(S_\ell)-1}) = (S_\ell(M''))(f(x_0),
\dotsc,f(x_{m(S_f)-1}$.
\endroster}
\ermn
For notational simplicity $S_\ell(M)(\bar b)$ is defined only for $\bar b$
with no repetitions, $m(S_i) \ge 0$ and $S_0,S_1,S_2,S_3$ are trivial.
\bn
\margintag{3.2}\ub{\stag{3.2} Minor Change}:  $x_{0,0}$ is connected to all $x_{\ell,i}$ for
$(\ell,i) \ne (0,0)$ (make the choice of representatives clear).
\bn
\margintag{3.3}\ub{\stag{3.3}}:  We demand $2^{h_2(i,\ell)} \le h_2(i,\ell+1)$ (when 
defined).
\bn
\margintag{3.4}\ub{\stag{3.4}}:  We define $M^2_{i,\ell(*)}$ and $M^2_i = M^2_{i,h_0(i)}$. 
\nl
The universe of $M^2_{i,\ell}$ is

$$
|M_i| \cup \bigcup \bigl\{ C_{i,\ell}:
\ell < \text{ Min}\{\ell(*),h_0(i)\} \bigr\} 
$$
\mn
where

$$
B_{i,\ell} = \bigl\{ \langle \ell,\bar b,a \rangle:\bar b \in {}^{m(S_\ell)}
|M_i| \text{ and } a \in |M_i| \bigr\}
$$
\mn
with no repetition for $u = \langle \ell,\bar b,a \rangle = \langle \ell(u),
\bar b(u),a(u) \rangle \in B_{i,\ell}$ let

$$
k(u) = h_1(\ell) (>0)
$$

$$
m[u] = \cases m(S_\ell) \times h_2(i,\ell) + (S_\ell(M_i))(\bar b)
&\text{ \ub{if} } a \notin \text{ Rang}(b[u]) \\
  m(S_\ell) \times h_2(i,\ell) + j &\text{ \ub{if} } a = b_j \endcases
$$
\mn
(note $j$ is unique as $\bar b$ with no repetition)

$$
t(u) = 2m(u),k(u).
$$
\mn
Lastly

$$
C_{i,\ell} = \bigl\{(u,c):u \in B_{i,\ell},c \in G_{k(u),t(u),m(u)} \bigr\}.
$$
\bn
\ub{Relations}:
\mr
\item "{$(a)$}"  those of $M^0_i$
\sn
\item "{$(b)-(d)$}"  as before.
\endroster
\bigskip

\definition{\stag{3.4A} Definition}  We add
\mr
\item "{$(e)$}"  $M^5_{i,\ell} = (M^2_i \restriction
(M_i \cup \{(u,c) \in M^2_i \backslash M_i:\ell(u) \le \ell\})$
$$
M^6_{i,\ell} = M^5_{i,\ell} + \Bbb N \restriction h(i,\ell +3).
$$
\endroster
\enddefinition
\bigskip

\definition{\stag{3.5} Definition}  $K_2 = \cup\{M^2_i / \cong:i < \omega\}$

$$
K^5_\ell = \cup\{M^5_{i,\ell} / \cong:i < \omega\}
$$

$$
K^0_\ell = \cup\{M^6_{i,\ell} / \cong:i < \omega\}.
$$
\enddefinition
\bigskip

\definition{\stag{3.6} Definition}  For $(K_0,K_2),(K_0,K^5_\ell),
(K^5_\ell,K^6_\ell)$ (continue as there).
\enddefinition
\bn
\par \noindent \llap{---$\!\!>$} MARTIN WARNS: Label 3.7 on next line is also used somewhere else (Perhaps should have used scite instead of stag?\par
\margintag{3.7}\ub{\stag{3.7}(1),(2)}:  The meaning is transformed; see \scite{1.7}(3).
\bigskip

\par \noindent \llap{---$\!\!>$} MARTIN WARNS: Label 3.7 on next line is also used somewhere else (Perhaps should have used scite instead of stag?\par
\proclaim{\stag{3.7}(3)Claim}  For each $\ell$ there are a LFP formulas
which interpret in $M^2_i$ a copy of $M^0_{i,\ell}$ say by a mapping $g_i$,
Dom$(g_i) = M^0_{i,\ell},g_i \restriction M^5_{i,\ell} =$ the identity
$g_i \restriction (\Bbb N \restriction h(i,\ell+3))$ has an image set of
equivalence classes (using elements of $M^6_{i,\ell+3} \backslash
M^6_{i,\ell+2}$.  In particular for some formula $\theta_\ell(x)$ for every
$i,(\forall b \in M^2_i)[d \in M^5_{i,\ell} \leftrightarrow M^2_i \models
\theta_\ell[d]]$.
\endproclaim
\bigskip

\demo{Proof}  Straight.
\enddemo
\bigskip

\proclaim{\stag{3.8A} Claim}  If $\varphi(\bar x) \in PFP(\tau(K_2))$, then
for some $\ell_\varphi < \omega$ and some formula $\psi_s(\bar x) \in PFP
(\tau(K^2_{\ell_\varphi}))(s < s^*)$ letting $\psi'_s(\bar x) = [\psi_s
(\bar x)]^{\theta_\ell(\omega)}$ (on $\theta_\ell(x)$ see above) 
(i.e. all free
variables and quantifications are demanded to satisfy $\theta_\ell(-)$ for
some Boolean combination $\psi(\bar x)$ of $\psi'_s(\bar x),s < s^*$ and
quantifier free formulas we have

$$
M \in K_2 \Rightarrow (\forall \bar x)(\bar \varphi(\bar x) \equiv
\psi(\bar x)).
$$
\endproclaim
\bigskip

\demo{Proof}  Like \scite{1.8} just easier.
\enddemo
\bigskip

\proclaim{\stag{3.8B} Claim}  Letting $n$ be fixed

$$
\align
{\Cal P}^n_i = \biggl\{ (Q,f):&Q \in {\Cal R}_i \text{ a } n
\text{-place relation on } M^5_{i,\ell},f \text{ a function with} \\
  &m(n) \text{-places (from the proof of \scite{1.8}, implicit) from} \\
  &M_i \text{ to } \Bbb N \restriction h_2(i,\ell +3) \text{ representing }
Q \biggr\}
\endalign
$$
\mn
give essential interpretations of $M^5_{i,\ell}$ in $M_i + \Bbb N \restriction
h_2(i,\ell +3)$ in the sense that:
\mr
\item "{$(*)_1$}"  for every first order $\varphi(Q,\bar x)$ variable
speaking on models in $K^5_\ell$ there is $\psi(f,\bar y,\bar z)$ first order
on speaking on $\{M + \Bbb N \restriction h_2(i,\ell +3;M\}$
\sn
\item "{$(*)_2$}"  if $(Q,f) \in {\Cal P}^n_i,Q^+ = \{\bar x:\varphi(Q,
\bar x)\}$ (relation) \nl
$f^+ = \{(\bar y,z):\psi(Q,\bar y,z)\}$ (graph of the function so
$\ell g(\bar y) = m_n$) \ub{then}
$(Q^+,f^+) \in {\Cal P}^{\ell g(\bar y)}_i$
\sn
\item "{$(*)_3$}"  for every $\psi(f,\bar y,z)$ above satisfying: if $f$ is
in $\{f:(\exists Q)(Qf) \in {\Cal P}^n_i\}$ then so is the function defined
by $\psi(f,\bar y,z)$ (i.e. it is syntactically clear not just semantically)
then there is $\varphi(Q,\bar x)$ such that $(*)_2$ above holds.
\endroster
\endproclaim
\bigskip

\proclaim{\stag{3.8} Claim}  Like \scite{1.8}.
\endproclaim
\bigskip

\proclaim{\stag{3.9} Claim}  Assume
\mr
\item "{$(*)$}"  if $\ell(*) < \omega,\psi = \psi(\bar x,y)$ a PFP formula
speaking on
$$
\{(M_i + \Bbb N \restriction h_2(i,\ell(*)+3),S_0(M_i),\dotsc,S_{\ell(*)}
(M_i)):i < \omega\}
$$
giving an $m$-place function from $M_i$ into $\Bbb N \restriction h_2(i,
\ell(*)+3)$, \ub{then} for some $\ell$ (generally $>> \ell(*)$), $S_\ell$ is
$\psi$.
\ermn
\ub{Then} in $K_2$, LFP and PFP are equivalent.
\endproclaim
\bigskip

\demo{Proof}  Just put together all.

So we have completed the proof of the theorem originally desired.
\enddemo

\newpage
    
REFERENCES.  
\bibliographystyle{lit-plain}
\bibliography{lista,listb,listx,listf,liste}

\enddocument

\bye